\documentclass[aap,nameyear,dvips]{arximspdf}
\usepackage{mathbh}
\usepackage{graphicx}

\doi{10.1214/07-AAP471}
\volume{18}
\issue{2}
\pubyear{2008}
\firstpage{664}
\lastpage{707}

\makeatletter
\newcommand{\eqref}[1]{\textup{(\ref{#1})}}
\def\nset{{\mathbb{N}}}
\def\rset{\mathbb R}

\def\rme{e}
\def\openset{\mathsf{O}}
\def\compactset{\mathsf{H}}
\def\unitsphere{\mathsf{S}}

\newcommand{\bepsilon}[2]{\bar{\epsilon}_{#2}}
\newcommand{\tepsilon}[2]{\tilde{\epsilon}_{#2}}
\newcommand{\teps}{\tilde{\epsilon}}
\newcommand{\pscal}[2]{\langle #1, #2 \rangle}
\newcommand{\1}{\mathbh{1}}
\def\Xset{\mathsf{X}} 
\def\Xsigma{\mathcal{X}}
\def\bphi{\bar{\Phi}} 
\def\bP{\bar{P}} 
\def\PE{\mathbb E}
\def\PP{\mathbb P}
\def\Q{\mathsf{Q}}
\newcommand{\eqdef}{\stackrel{\mathrm{def}}{=}}
\newcommand{\tvnorm}[1]{ \|#1 \|_{\mathrm{TV}}}
\newcommand{\fnorm}[2]{\|#1 \|_{#2}}
\newcommand{\Hphi}{H_\phi}
\def\lleb{\lambda^{\mathrm{Leb}}}

\newtheorem{theo}{Theorem}[section]
\newtheorem{lem}[theo]{Lemma}
\newtheorem{coro}[theo]{Corollary}
\newtheorem{prop}[theo]{Proposition}
\newproclaim{defi}[theo]{Definition}
\newproclaim{rem}{Remark}
\newproclaim{example}{Example}
\newcommand{\QQq}[2]{\mathbb{Q}^{#1}_{#2}}
\newcommand{\continuousf}[1]{\mathsf{C}(#1)}
\makeatother

\begin{document}
\begin{frontmatter}

\title{The ODE method for stability   of skip-free Markov chains with applications to MCMC\thanksref{T1}}
\runtitle{The ODE Method for Markov Chain Stability}
\thankstext{T1}{Supported in part by the National Research Agency (ANR) under the program ``ANR-05-BLAN-0299.''}

\begin{aug}
\author[A]{\fnms{Gersende} \snm{Fort}\ead[label=e1]{gfort@tsi.enst.fr}\corref{}},
\author[B]{\fnms{Sean} \snm{Meyn}\ead[label=e3]{meyn@uiuc.edu}},
\author[A]{\fnms{Eric} \snm{Moulines}\ead[label=e2]{moulines@tsi.enst.fr}} \and
\author[C]{\fnms{Pierre} \snm{Priouret}\ead[label=e4]{priouret@ccr.jussieu.fr}}
\runauthor{G. Fort, S. Meyn, E. Moulines and P. Priouret}
\affiliation{T\'{e}l\'{e}com Paris\textup{,} CNRS\textup{,}
University of Illinois at Urbana-Champaign,\\
T\'{e}l\'{e}com Paris\textup{,} CNRS
and Universit\'e Pierre et Marie Curie}
\address[A]{G. Fort\\
E. Moulines\\
Laboratoire Traitement\\
\quad et Communication de l'Information\\
T\'{e}l\'{e}com Paris/CNRS\\
46 rue Barrault\\
75634 Paris C\'{e}dex 13\\
France\\
\printead{e1}\\
\phantom{\textsc{E-mail}: } \printead*{e2}}
\address[B]{S. Meyn\\
Department of Electrical\\
\quad and Computer Engineering\\
and\\
the Coordinated Sciences Laboratory\\
University of Illinois at Urbana-Champaign\\
Urbana, Illinois 61801\\
USA\\
\printead{e3}}
\address[C]{P. Priouret\\
Laboratoire de Probabilit\'{e}s\\
\quad et Mod\`{e}les  Al\'{e}atoires\\
Universit\'{e} Pierre et Marie Curie\\
Bo\^{\i}te courrier 188\\
75252 Paris Cedex 05\\
France\\
\printead{e4}}
\end{aug}

\received{\smonth{6} \syear{2006}}
\revised{\smonth{5} \syear{2007}}

\begin{abstract}
Fluid limit techniques have become a central tool to analyze queueing
networks over the last decade, with applications to performance
analysis, simulation and optimization.

In this paper, some of these techniques are extended to a general class
of skip-free Markov chains.  As in the case of queueing models, a fluid
approximation is obtained by scaling time, space and the initial
condition by a large constant.  The resulting fluid limit is the
solution of an ordinary differential equation (ODE) in ``most'' of the
state space.  Stability and finer ergodic properties for the stochastic
model then follow from stability of the set of fluid limits.  Moreover,
similarly to the queueing context where fluid models are routinely used
to design control policies, the structure of the limiting ODE in this
general setting provides an understanding of the dynamics of the Markov
chain.  These results are illustrated through application to Markov
chain Monte Carlo  methods.
\end{abstract}


\end{frontmatter}

The use of ordinary differential equations (ODE) to analyze
Markov chains was first suggested by \citet{kurtz1970}.
This idea was later refined by \citet{newell1982}, who introduced the so-called
\textit{fluid approximations} with applications to queueing networks.
Since the 1990s, fluid models have been
used to address delay in complex networks [\citet{cruz1991}] and bottleneck
analysis [\citet{chenmandelbaum1991}].  The latter work followed an already
extensive research program on diffusion approximations for networks [see
\citet{harrison2000}, \citet{whitt2002}, \citet{chenyao2001} and the references therein].

The purpose of this paper is to extend fluid limit techniques to a general
class of discrete-time Markov chains $\{ \Phi_k \}$ on a $d$-dimensional
Euclidean state space $\Xset$.  Recall that a Markov chain is called
\emph{skip-free} if the increments $( \Phi_{k+1}-\Phi_k)$ are uniformly bounded
in norm by a deterministic constant for each $k$ and each initial condition.
For example, Markov chain models of queueing systems are typically skip-free.
Here, we consider a relaxation of this assumption in which the increments are
assumed to be bounded in an $L^p $-sense.  Consequently, we find that the chain can
be represented by the additive noise model
\begin{equation} \label{eq:random-walk-type}
\Phi_{k+1} = \Phi_k + \Delta(\Phi_k)+ \epsilon_{k+1},
\end{equation}
where $\{\epsilon_k\}$ is a martingale increment sequence w.r.t. the natural
filtration of the process $\{\Phi_k\}$ and $\Delta\colon\Xset \to\Xset$ is
bounded.  Associated with this chain, we consider the sequence of
continuous-time processes
\begin{eqnarray}  \label{eq:InterpolatedProcessIntro}
  \eta_{r}^{\alpha}({t};{x})  \eqdef   r^{-1} \Phi_{\lfloor t r^{1+ \alpha} \rfloor},\qquad
  \eta_{r}^{\alpha}({t};{0})= r^{-1} \Phi_0 = x,
  \nonumber
  \\[-8pt]
  \\[-8pt]
\eqntext{ r \geq 0, \alpha \geq 0, x \in \Xset,}
\end{eqnarray}
obtained by interpolating and scaling the Markov chain in space and time.  A
fluid limit is obtained as a subsequential weak limit of a sequence $\{
\eta_{r_n}^{\alpha}(\cdot;x_n)\}$, where $\{r_n\}$ and $\{x_n\}$ are two
sequences such that $\lim_{n \to \infty} r_n = \infty$ and $\lim_{n \to \infty}
x_n =x$.  The set of all such limits is called the \textit{fluid limit model}.
In queueing network applications, a fluid limit is easy to interpret in terms
of mean flows; in most situations, it is a solution of a deterministic set of
equations depending on network characteristics as well as the control policy
[see, e.g., \citet{chenmandelbaum1991}, \citet{dai1995}, \citet{daimeyn1995}, \citet{chenyao2001},
\citet{meyn-book2006}]. The existence of limits and the continuity of the fluid limit model may be
established under general conditions on the increments (see
Theorem~\ref{theo:weak-compacity}).

The fact that stability of the fluid limit model implies stability of the
stochastic network was established in a limited setting in
\citet{malyshevmenshikov1979}. This was extended to a very broad class of multiclass
networks by \citet{dai1995}.  A key step in the proof of these results is a
multi-step state-dependent version of Foster's criterion introduced in
\citet{malyshevmenshikov1979} for countable state space models, later
extended to general state space models in Meyn and Tweedie (\citeyear{meyntweedie1993,meyntweedie1994d}).
The main result of \citet{dai1995} only established positive recurrence.
Moments and rates of convergence to stationarity of the Markovian network model
were obtained in \citet{daimeyn1995}, based on an extension of
\citet{meyntweedie1994d} using the subgeometric $f$-ergodic theorem in
\citet{tuominentweedie1994} [recently extended and simplified in work of
\citet{doucfortmoulinessoulier2004}].  Converse theorems have appeared in
\citet{daiweiss1996}, \citet{dai1996}, \citet{meyn1995} that show that, under rather
strong conditions, instability of the fluid model implies transience of the stochastic
network.  The counterexamples in
\citet{gamarnikhasenbein2005}, \citet{daihasenbeinvandevate2004} show that some
additional conditions are necessary to obtain a converse.

Under general conditions, including the generalized skip-free assumption, a
fluid limit $\eta$ is a \emph{weak} solution (in a sense given below)
to the homogeneous ODE
\begin{equation}
\label{eq:definition-ODE}
\dot{\mu} = h(\mu).
\end{equation}
The vector field  $h$ is defined as a radial limit of the function $\Delta$ appearing in \eqref{eq:random-walk-type} under appropriate renormalization.

Provided that the increments $\{ \epsilon_k \}$ in the decomposition
\eqref{eq:random-walk-type} are tight in $L^p$, stability of the fluid limit
model implies finite moments in steady state, as well as polynomial rates of
convergence to stationarity; see Theorem \ref {theo:(f,r)-ergodicity}.

One advantage of the ODE approach over the usual Foster--Lyapunov approach to
stability is that the ODE model provides insight into Markov chain dynamics.
In the queueing context, the ODE model has many other
applications, such as simulation variance reduction
[\citet{hendersonmeyntadic2003}] and optimization [\citet{chenmeyn1999}].

The remainder of the paper is organized as follows. Section \ref{subsec:Fluid
  Limit: definitions} contains notation and assumptions, along with a
construction of the fluid limit model.  The main result is contained in Section
\ref{subsec:Stability of Fluid Limits and Markov Chain Stability}, where it is
shown that stability of the fluid limit model implies the existence of
polynomial moments as well as polynomial rates of convergence to stationarity
[known as $(f,r)$-\textit{ergodicity}].

Fluid limits are characterized in Section \ref{subsec:Characterization of the
  fluid limits}. Proposition~\ref{prop:characterization-fluid-limit} provides
conditions that guarantee that a fluid limit coincides with the weak solutions
of the ODE \eqref{eq:definition-ODE}.

These results are applied to establish $(f,r)$-ergodicity of the random walk
Metropolis--Hastings algorithm for superexponential densities in
Section~\ref{subsec:Super-exponential target densities} and subexponential
densities in Section~\ref{sec:weibulian-tails}.  In Examples
\ref{example:Mixture of Gaussian Densities} and
\ref{example:mixture-weibulian}, the fluid limit model is stable and any fluid
limit is a weak solution of the ODE \eqref{eq:definition-ODE}, yet some fluid
limits are nondeterministic.

The conclusions contain proposed extensions, including diffusion limits of the
form obtained in \citet{harrison2000}, \citet{whitt2002}, \citet{chenyao2001} and application
of ODE methods for variance reduction in simulation and MCMC.

\section{Assumptions and statement of the results}
\label{sec:theory-fluid-limit}

\subsection{Fluid limit\textup{:} definitions}
\label{subsec:Fluid Limit: definitions}
We consider a Markov chain $\bolds{\Phi} \eqdef \{ \Phi_k \}_{k \geq 0}$
on a $d$-dimensional Euclidean space $\Xset$ equipped with its Borel
sigma-field $\Xsigma$.  We denote by $\{\mathcal{F}_k \}_{k \geq 0}$ the
natural filtration. The distribution of $\bolds{\Phi}$ is specified by its
initial state $\Phi_0=x \in \Xset$ and its transition kernel $P$. We write
$\PP_x$ for the distribution of the chain conditional on the initial state
$\Phi_0=x$ and $\PE_x$ for the corresponding expectation.

Denote by $\mathsf{C}(\rset^+,\Xset)$ the space of continuous
$\Xset$-valued functions on the infinite time interval $[0, \infty)$.  We equip
$\mathsf{C}(\rset^+,\Xset)$ with the local uniform topology.  Denote
by $\mathsf{D}(\rset^+,\Xset)$ the space of $\Xset$-valued
right-continuous functions with left limits on the infinite time interval $[0,
\infty)$, hereafter \textit{c\`adl\`ag functions}.  This space is endowed with the
Skorokhod topology.  For $0<T<+\infty$, denote by
$\mathsf{C}([0,T],\Xset)$ (resp.  $\mathsf{D}([0,T],\Xset)$)
the space of $\Xset$-valued continuous functions (resp.  c\`adl\`ag functions)
defined on $[0,T]$, equipped with the uniform (resp. Skorokhod) topology.

For $x \in \Xset$, $\alpha \geq 0$ and $r > 0$,  consider the interpolated process
\begin{equation}
\label{eq:definition:fluid:limit}
\eta_r^{\alpha}({t};{x})  \eqdef   r^{-1} \Phi_{\lfloor t r^{1+ \alpha} \rfloor},\qquad  \eta_r^{\alpha}({t};{0})= r^{-1} \Phi_0 = x,
\end{equation}
where  $\lfloor \cdot \rfloor$ stands for the lower integer part.
Denote by $\mathbb{Q}^{\alpha}_{r;x}$ the image probability on
$\mathsf{D}(\rset^+,\Xset)$ of $\PP_x$ by $\eta_r^{\alpha}(\cdot;x)$. In words, the renormalized process is obtained by
scaling the Markov chain in space, time and initial condition.
This is made precise in the following definition.

\begin{defi}[($\alpha$-\textit{fluid limit})]
  Let $\alpha \geq 0$ and $x \in \Xset$.  A probability measure
  $\mathbb{Q}^{\alpha}_{x}$ on $\mathsf{D}(\rset^+,\Xset)$ is said to be an
  \emph{$\alpha$-fluid limit} if there exist sequences of scaling factors $\{
  r_n \} \subset \rset_+$ and initial states $\{x_n \} \subset \Xset$
  satisfying $\lim_{n \to \infty} r_n = +\infty$ and $\lim_{n \to \infty} x_n=
  x$ such that $\{ \mathbb{Q}_{r_n;x_n}^{\alpha}\}$ converges weakly to
  $\mathbb{Q}^{\alpha}_{x}$ on $\mathsf{D}(\rset^+,\Xset)$ (denoted
  $\mathbb{Q}_{r_n;x_n}^{\alpha} \Rightarrow \mathbb{Q}^{\alpha}_{x}$).
\end{defi}

The set $\{ \mathbb{Q}^{\alpha}_{x}, x \in \Xset\}$ of all such limits is referred to as
the $\alpha$-\emph{fluid limit model}. An $\alpha$-fluid limit $\mathbb{Q}^{\alpha}_{x}$
is said to be \emph{deterministic} if there exists a function $g \in
\mathsf{D}(\rset^+,\Xset)$ such that $\mathbb{Q}^{\alpha}_{x} = \delta_{g}$, the
Dirac mass at $g$.

Assume that $\PE_x[ |\Phi_1| ] < \infty$ for all $x \in \Xset$, where
$|\cdot|$ denotes the Euclidean norm, and consider the decomposition
\begin{equation}
\label{eq:decomposition-Phi}
\Phi_k = \Phi_{k-1} + \Delta(\Phi_{k-1}) + \epsilon_k,   \qquad  k \geq 1,
\end{equation}
where
\begin{eqnarray}
\label{eq:DefDelta}
\Delta(x) &\eqdef& \PE_x[\Phi_1 - \Phi_0]= \PE_x  [ \Phi_1  ]- x\qquad  \mbox{for all $x \in \Xset$},
\\
\label{eq:definition-epsilon}
\epsilon_k &\eqdef& \Phi_k - \PE [\Phi_k \vert \mathcal{F}_{k-1}]\hspace*{23.2mm}  \mbox{for all $k \geq 1$}.
\end{eqnarray}
In the sequel, we assume the following.
\begin{enumerate}[B2.]
\item[B1.] There exists $p>1$ such that $\lim_{K \to \infty} \sup_{x \in \Xset} \PE_x[ |\epsilon_1|^p \1 \{ |\epsilon_1| \geq K \}] = 0$.

\item[B2.]  There exists $\beta \in [0, 1 \wedge (p -1)  )$ such that
$N(\beta,\Delta) \eqdef \sup_{x \in \Xset}  \{ (1 + |x|^\beta)\times | \Delta(x)|  \}  < \infty$.
\end{enumerate}

\begin{theo}\label{theo:weak-compacity}
Assume \textup{B1} and \textup{B2}.  Then, for all $0 \leq \alpha \leq \beta$
and any sequences $\{ r_n \} \subset \rset_+$ and $\{x_n \} \subset \Xset$ such
that $\lim_{n \to \infty} r_n = +\infty$ and $\lim_{n \to \infty} x_n = x$,
there exists a probability measure $\mathbb{Q}^{\alpha}_{x}$ on
$\mathsf{C}(\rset_+,\Xset)$ and subsequences $\{ r_{n_j} \} \subseteq
\{r_n \}$ and $\{ x_{n_j} \} \subseteq \{x_n\}$ such that
$\mathbb{Q}_{r_{n_j};x_{n_j}}^\alpha  \Rightarrow \mathbb{Q}^{\alpha}_{x}$.  Furthermore, for
all $ 0 \leq \alpha < \beta$, the $\alpha$-fluid limits are \emph{trivial} in
the sense that $\mathbb{Q}^{\alpha}_{x}= \delta_g$ with $g(t) \equiv x$.
\end{theo}

Note that for any $x \in \Xset$ and $0 \leq \alpha \leq \beta$,
we have $\mathbb{Q}^{\alpha}_{x}( \eta, \eta(0)= x) = 1$, showing that $x$ is the initial point of the fluid limit.

\subsection{Stability of fluid limits and Markov chain stability}
\label{subsec:Stability of Fluid Limits and Markov Chain Stability}
There are several notions of stability that have appeared in the literature [see
\citet{meyn2001}, Theorem~3] and the surrounding discussion. We adopt the
notion of stability introduced in \citet{stolyar1995}.

\begin{defi}[(\textit{Stability})] The $\alpha$-fluid limit model  is said to be \emph{stable} if there exist $T > 0$
  and $\rho < 1$ such that for any $x \in \Xset$ with $|x|=1$,
\begin{equation}
\label{eq:stability:fluidlimit-small}
\mathbb{Q}^{\alpha}_{x}   \biggl( \eta \in \mathsf{D}(\rset_+,\Xset),   \inf_{0 \leq t \leq T} |\eta(t)| \leq \rho  \biggr) = 1.
\end{equation}
\end{defi}

Let $f \dvtx  \Xset \to [1,\infty)$ and $L_\infty^f$ denote the vector space of all
measurable functions $g$ on $\Xset$ such that $\sup_{x \in \Xset} |g(x)|/f(x)$
is finite. $L_\infty^f$ equipped with the norm $| g |_f \eqdef \sup_{x \in
  \Xset} |g(x)|/f(x)$ is a Banach space.

Denote by $\fnorm{\cdot}{f}$ the $f$-total variation norm, defined for any
finite signed measure $\nu$ as
$\fnorm{\nu}{f} = \sup_{|g| \leq f}|\nu(g)|$.

We recall some basic definitions related to Markov chains on general state space; see
\citet{meyntweedie1993} for an in-depth presentation. A~chain is said
to be \textit{phi-irreducible} if there exists a $\sigma$-finite measure $\phi$ such that
$\sum_{n \geq 0} P^n(x,A) > 0$ for all $x \in \Xset$ whenever $\phi(A) > 0$.  A
set $C \in \Xsigma$ is $\nu_m$-\textit{small} if there exist a nontrivial measure
$\nu_m$ and a positive integer $m$ such that such that $P^m(x, \cdot) \geq
\1_C(x) \nu_m(\cdot)$.  \textit{Petite} sets are a generalization of small sets: a set
$C$ is said to be petite if there exists a distribution $a$ on the positive integers
and a distribution $\nu$ such that $\sum_{n \geq 0} a(n) P^n(x,\cdot) \geq
\1_C(x) \nu(\cdot)$.  Finally, an \textit{aperiodic chain} is a chain such that the
greatest common divisor of the set
\[
\{m, C   \mbox{is $\nu_m$-small and $\nu_m = \delta_m \nu$ for some $\delta_m
>0$} \},
\]
is one, for some small set $C$.  For a phi-irreducible aperiodic chain, the petite sets
are small [\citet{meyntweedie1993}, Proposition~5.5.7].

Let $\{ r(n) \}_{n \in \nset}$ be a sequence of positive real numbers. An
aperiodic phi-irreducible positive Harris chain with stationary distribution
$\pi$ is called\break \mbox{$(f,r)$-}\textit{ergodic} if
\[
\lim_{n \to \infty} r(n) \fnorm{P^n(x,\cdot) - \pi}{f} = 0
\]
for all $x \in \Xset$. If $P$ is positive Harris recurrent with invariant
probability $\pi$, then the fundamental kernel $Z$ is defined as $Z \eqdef
(\mathrm{Id} - P + \Pi)^{-1}$, where the kernel $\Pi$ is $\Pi(x,\cdot) \equiv
\pi(\cdot)$ for all $x \in \Xset$ and $\mathrm{Id}$ is the identity kernel.
For any measurable function $g$ on $\Xset$, the function $\hat{g}= Z g$ is a
solution to the Poisson equation, whenever the inverse is well defined [see
\citet{meyntweedie1993}].

The following theorem may be seen as an extension of [\citet{daimeyn1995},
Theorem~5.5], which relates the stability of the fluid limit to the
$(f,r)$-ergodicity of the original chain.

\begin{theo}
\label{theo:(f,r)-ergodicity}
Let $\{\Phi_k \}_{k \in \nset}$ be a phi-irreducible and aperiodic Markov chain
such that compact sets are petite. Assume
\textup{B1}~and~\textup{B2}
and that the $\beta$-fluid limit model is stable. Then, for any $1 \leq q \leq
(1+\beta)^{-1} p$,
\begin{longlist}[(ii)]
\item[(i)] the Markov chain $\{ \Phi_k \}_{k \in \nset}$ is
  $ (f^{(q)}, r^{(q)}  )$-ergodic with $f^{(q)}(x) \eqdef 1+|x|^{p-q(1+\beta)}$ and
  $r^{(q)}(n)= n^{q-1}$;

\item[(ii)] the fundamental kernel $Z$ is a bounded linear transformation from
  $L_\infty^{f^{(q)}}$ to $L_\infty^{f^{(q-1)}}$.
\end{longlist}
\end{theo}

\subsection{Characterization of the fluid limits}
\label{subsec:Characterization of the fluid limits}
Theorem \ref{theo:(f,r)-ergodicity} relates the ergodicity of the Markov chain
to the stability of the fluid limit and raises the question: \emph{how can we
  determine if the  $\beta$-fluid model is stable}? To answer this
question, we first characterize the set of fluid limits.

In addition to assumptions B1--B2, we require   conditions on the limiting behavior of the
function $\Delta$.
\begin{enumerate}[B3.]
\item[B3.] There exist an open cone  $\openset
\subseteq \Xset \setminus \{0\}$ and  a continuous function $\Delta_\infty: \openset \to \Xset$
such that, for any compact subset $\compactset \subseteq \openset$,
\[
\lim_{r \to +\infty} \sup_{x \in \compactset}    \big| r^{\beta} |x|^{\beta} \Delta(rx) - \Delta_\infty(x)  \big| = 0,
\]
where $\beta$ is given by B2.
\end{enumerate}

The easy
situation is when \mbox{$\openset = \Xset \setminus \{0\}$,} in which case the radial
limit\break $ \lim_{r \to \infty} r^\beta |x|^\beta \Delta(rx)$ exists for $x \ne 0$.
Though this condition is met in examples of interest, there are several
situations for which the radial limits do not exist for directions belonging to
some low-dimensional manifolds of the unit sphere. Let $h$ be given by
\begin{equation}
\label{eq:DefH}
h(x) \eqdef |x|^{-\beta}   \Delta_\infty (x).
\end{equation}
A function $\mu\dvtx  I \to \Xset$ (where $I \subset \rset^+$ is an interval
which can be open or closed, bounded or unbounded) is said to be a
\emph{solution of the ODE \eqref{eq:definition-ODE}} on $I$ with initial
condition $x$ if $\mu$ is continuously differentiable on $I$ for all $t \in
I$, $\mu (t) \in \openset$, $\mu(0) =x$ and $\dot{\mu}(t)= h \circ
\mu(t)$. The following theorem shows that the fluid limits restricted to $\openset$ evolve deterministically and, more precisely,
that their supports on $\openset$ belong to the flow of the ODE.

\begin{prop}\label{prop:characterization-fluid-limit}
Assume \textup{B1, B2} and \textup{B3}. For any $ 0 \leq s \leq t $, define
\begin{equation}
\mathsf{A}(s,t) \eqdef  \{ \eta \in \mathsf{C}(\rset^+,\Xset) \dvtx \eta(u) \in \openset      \mbox{ for all } u \in [s,t]  \}.
\end{equation}
Then, for any $x \in \Xset$ and any $\beta$-fluid limit~$\QQq{\beta}{x}$, on $\mathsf{A}(s,t)$,
\[
\sup_{s \leq u \leq t}  \bigg| \eta(u) - \eta(s) - \int_s^u h \circ \eta(v)\,dv  \bigg|=
0,\qquad \mathbb{Q}^\beta_x\mbox{-a.s.}
\]
\end{prop}

Under very weak additional conditions, one may assume that the
solutions of the ODE \eqref{eq:definition-ODE} with initial condition $x \in
\openset$ exist and are unique on a nonvanishing interval $[0,T_{x}]$.  In
such a case, Proposition~\ref{prop:characterization-fluid-limit} provides a handy
description of the fluid limit.
\begin{enumerate}[B4.]
\item[B4.] Assume that for all $x \in \openset$, there exists
  $T_{x} > 0$ such that the ODE \eqref{eq:definition-ODE} with initial
  condition $x$ has a \emph{unique} solution, denoted
  $\mu(\cdot;x)$ on an interval $[0,T_{x}]$.
\end{enumerate}

Assumption B4 is satisfied if
$\Delta_\infty$ is locally Lipschitz on  $\openset$; in such a case,
$h$ is locally Lipschitz on $\openset$ and it then follows from
classical results on the existence of solutions of the ODE [see,  e.g.,
\citet{verhulst1996}] that for any $x \in \openset$, there exists
$T_x > 0$ such that, on the interval $[0,T_x]$, the ODE
\eqref{eq:definition-ODE} has a unique solution $\mu$
with initial condition $\mu(0)= x$.  In addition, if the
ODE \eqref{eq:definition-ODE} has  two solutions,
$\mu_1$~and~$\mu_2$, on an interval $I$
which satisfy $\mu_1(t_0)=  \mu_2(t_0)
= x_0$ for some $t_0 \in I$, then  $\mu_1(t)=
\mu_2(t)$ for any $t \in I$.

An elementary application of Proposition~\ref{prop:characterization-fluid-limit}
shows that under this additional assumption, a fluid limit starting at $x_0
\in \openset$ coincides with the solution of the ODE \eqref{eq:definition-ODE}
with initial condition $x_0$ on a nonvanishing interval.

\begin{theo}
\label{theo:fluidlimit-Lipsch-NeighLimit}
Assume \textup{B1--B4}. Let
$x \in \openset$. There then exists $T_{x} > 0$ such that $\mathbb{Q}^\beta_x =
\delta_{\mu(\cdot;x)}$ on $\mathsf{D}([0,T_{x}],\Xset)$.
\end{theo}

As a corollary of Theorem~\ref{theo:fluidlimit-Lipsch-NeighLimit}, we
have the following.

\begin{coro}
  Assume that $\openset= \Xset \setminus \{0\}$ in B3.
  Then all $\beta$-fluid limits are deterministic and solve the ODE
  \eqref{eq:definition-ODE}.  Furthermore, for any $\epsilon>0$ and $x \in
  \Xset$, and any sequences $\{r_n\} \subset \rset_+$ and $\{x_n\} \subset
  \Xset$ such that $\lim_{n \to \infty} r_n = + \infty$ and $\lim_{n \to
    \infty} x_n = x$,
\[
\lim_n \PP_{r_n x_n}  \biggl( \sup_{0 \leq t \leq T_{x}}  | \eta_{r_n}^{\beta}({t};{x_n}) - \mu(t;x)  | \geq \epsilon  \biggr) =0.
\]
\end{coro}

Hence, the fluid limit depends only on the initial value $x$ and does not
depend upon the choice of the sequences $\{r_n\}$ and $\{x_n\}$.

The last step is to relate the stability of the fluid limit [see
(\ref{eq:stability:fluidlimit-small})] to the behavior of the solutions of the
ODE, when such solutions are well defined. From the discussion above, we may
deduce a first elementary stability condition. Assume that
B3 holds with $\openset = \Xset \setminus \{0\}$.  In
this case, the fluid limit model is stable if there exist $\rho < 1$ and $T <
\infty$ such that, for any $|x|=1$, $\inf_{[0,T]} |\mu(\cdot;x)| <
\rho$,  that is, the solutions of the ODE enter a sphere of radius $\rho < 1$ before
a given time $T$.

\begin{theo}
\label{theo:stability-fluid-limit-deterministic-smooth}
Let $\{\Phi_k \}_{k \in \nset}$ be a phi-irreducible and aperiodic Markov chain such that compact sets are petite.
Let $\rho$, $0 < \rho < 1$ and $T > 0$.  Assume that
\textup{B1--B4}
hold with $\openset= \Xset \setminus \{0\}$.  Assume, in addition,
that for any $x$ satisfying $|x|=1$, the solution
$\mu(\cdot;x)$ is such that $\inf_{[0,T \wedge T_{x}]} |\mu(\cdot;x)| \leq \rho$.  Then, the
$\beta$-fluid limit model is stable and the conclusions of Theorem
\textup{\ref{theo:(f,r)-ergodicity}} hold.
\end{theo}

When B3 holds for a strict subset of the
state space $\openset \subsetneq \Xset \setminus\{0\}$, the
situation is more difficult because some fluid limits are not
solutions of the ODE.  Regardless, under general assumptions,
stability of the ODE implies stability of the fluid limit model.

\begin{theo}
\label{theo:stability-fluid-limit-deterministic-nonsmooth}
Let $\{\Phi_k \}_{k \in \nset}$ be a phi-irreducible and aperiodic Markov chain such that compact sets are petite.
Assume that \textup{B1--B4}
hold with $\openset \subsetneq \Xset \setminus \{0\}$. Assume, in addition,
that:
\begin{longlist}[(iii)]
\item[(i)] \label{item:IntersectCGamma}
  there exists $T_0 > 0$ such that for any $x$, $|x|=1$, and for any $\beta$-fluid limit~$\mathbb{Q}^\beta_x$,
  \begin{equation}    \label{eq:IntersectCGamma}
     \mathbb{Q}^\beta_x \bigl(\eta\dvtx \eta([0,T_0]) \cap \openset \neq \varnothing  \bigr) =1;
  \end{equation}

\item[(ii)] \label{item:definition:TK} for any $K>0$, there exist $T_K > 0$ and $0 < \rho_K <1$ such that for any $x \in \openset$, $|x| \leq K$,
  \begin{equation}
    \label{eq:DefiTK}
    \inf_{[0, T_K \wedge T_x]}| \mu(\cdot;x) | \leq \rho_K;
  \end{equation}

\item[(iii)] \label{item:definition:OmegaC} for any compact set $\compactset \subset
  \openset$ and any $K$,
  \[
  \Omega_{\compactset} \eqdef   \{
    \mu(t;x)\dvtx  x \in \compactset, t \in [0,T_x \wedge T_K]  \}
\]
is a compact subset of $\openset$.
\end{longlist}
Then, the  $\beta$-fluid model is stable and the conclusions of
Theorem \textup{\ref{theo:(f,r)-ergodicity}} hold.
\end{theo}

Condition (i) implies that each
$\beta$-fluid limit reaches the set $\openset$ in a finite time. When the
initial condition $x \ne 0$ does belongs to $\openset$, this condition is
automatically fulfilled. When $x$ does not belong to
$\openset$, this condition typically requires that there is a force driving
the chain into $\openset$. The verification of this property
generally requires some problem-dependent and sometimes intricate constructions
(see,  e.g.,  Example~\ref{example:Mixture of Gaussian Densities}).
Condition (ii) implies that the solution
$\mu(\cdot;x)$ of the ODE with initial point $x \in \openset$ reaches
a ball inside the unit sphere before approaching the singularity. This also
means that the singular set is repulsive for the solution of the ODE.

\section{The ODE method for the Metropolis--Hastings algorithm}
\label{sec:The ODE method for the Metropolis-Hastings algorithm}
The Metro\-polis--Hastings (MH) algorithm [see \citet{robertcasella1999} and the
references therein] is a popular computational method for generating samples
from virtually any distribution $\pi$. In particular, there is no need for the
normalizing constant to be known and the space $\Xset = \rset^d$ (for some
integer $d$) on which it is defined can be high-dimensional.  The method
consists of simulating an ergodic Markov chain $\{\Phi_k\}_{k \geq 0}$ on
$\Xset$ with transition probability $P$ such that $\pi$ is the stationary
distribution for this chain, that is, $\pi P = \pi$.

The MH algorithm requires the choice of a \emph{proposal kernel} $q$. In order
to simplify the discussion, we will here assume that $\pi$ and $q$ admit
densities with respect to the Lebesgue measure $\lleb$, denoted (with an abuse
of notation) $\pi$ and $q$ hereafter. We denote by $\Q$ the probability defined
by $\Q(A)= \int_A q(y) \lleb(dy)$. The role of the kernel $q$ consists of
proposing potential transitions for the Markov chain $\{ \Phi_k \}$. Given that
the chain is currently at $x$, a candidate $y$ is accepted with probability
$\alpha (x,y)$, defined as
$\alpha (x,y)= 1\wedge \frac{\pi (y)}{\pi (x)}\frac{q(y,x)}{q(x,y)}$.
Otherwise it is rejected and the Markov
chain stays at its current location $x$.  The transition kernel $P$ of this
Markov chain takes the form, for $x \in \Xset$ and $A \in
\mathcal{B}(\Xset)$,
\begin{eqnarray} \label{eq:kernelSRWM}
P(x,A)&=&  \int_{A-x} \alpha(x,x+y)q(x,x+y) \lleb(dy)
\nonumber
\\[-8pt]
\\[-8pt]
\nonumber
&&{} + \1_A(x) \int_{\Xset-x} \{ 1-\alpha(x,x+y) \} q(x,x+y) \lleb(dy),
\end{eqnarray}
where $A - x \eqdef \{ y \in \Xset, x+y \in A\}$.
The Markov chain $P$ is reversible with respect to $\pi$ and therefore admits
$\pi$ as invariant distribution. For the purpose of illustration, we focus on
the symmetric increments random walk MH algorithm (hereafter SRWM), in which
$q(x,y)=q(y-x)$ for some symmetric distribution $q$ on $\Xset$.  Under
these assumptions, the acceptance probability simplifies to $\alpha(x,y) = 1
\wedge [\pi(y)/\pi(x)]$.  For any measurable function $W \dvtx  \Xset \to \Xset$,
\begin{eqnarray*}
\PE_x  [ W(\Phi_1)  ] - W(x) &=&  \int_{\mathsf{A}_x} \{ W(x+y) - W(x) \} q(y) \lleb(dy)
\\
&&{} + \int_{\mathsf{R}_x} \{ W(x+y) - W(x) \} \frac{\pi(x+y)}{\pi(x)} q(y) \lleb(dy),
\end{eqnarray*}
where $\mathsf{A}_x \eqdef \{y \in \Xset, \pi(x+y) \geq \pi(x) \}$ is the acceptance
region (moves toward $x+\mathsf{A}_x$ are accepted with probability one) and $\mathsf{R}_x
\eqdef \Xset \setminus \mathsf{A}_x$ is the potential rejection region.
From \citet{robertstweedie1996}, Theorem~2.2,  we obtain the following basic result.

\begin{theo}
\label{theo:PetiteIrred}
Suppose that the target density $\pi$ is positive and continuous and that $q$
is bounded away from zero,  that is, there exist $\delta_q> 0$ and $\epsilon_q> 0$
such that $q(x) \geq \epsilon_q$ for $|x| \leq \delta_q$. Then, the
random-walk-based Metropolis algorithm on $\{ \Xset, \Xsigma\}$
is $\lleb$-irreducible, aperiodic and every nonempty bounded set is small.
\end{theo}

In the sequel, we assume that $q$ has a moment of order $p > 1$.
To apply the results presented in Section~\ref{sec:theory-fluid-limit}, we must
first compute $\Delta(x)= \PE_x[ \Phi_1] - x$,  that is, to set $W(x)=x$ in the
previous formula. Since $q$ is symmetric and therefore zero-mean, the previous
reduces to
\begin{equation}
  \label{eq:DefinitionDelta}
  \Delta(x) = \int_{\mathsf{R}_x} y    \biggl( \frac{\pi(x+y)}{\pi(x)} -1  \biggr) q(y) \lleb(dy).
\end{equation}
Note that, for any $x \in \Xset$, $|\epsilon_1| \leq |\Phi_1-\Phi_0| + m \PP_x$-a.s., where $m = \int |y| q(y)\times \lleb(dy)$. Therefore, for any $K > 0$,
\begin{eqnarray*}
  \PE_x  [|\epsilon_1|^p \1 \{ |\epsilon_1| \geq K \}  ] &\leq&
  2^p \PE_x[ (|\Phi_1 - \Phi_0|^p + m^p) \1 \{ |\Phi_1-\Phi_0| \geq K-m\}]
  \\
  &\leq& 2^p \int |y|^p \1 \{ |y| \geq K-m \} q(y) \lleb(dy),
\end{eqnarray*}
showing that assumption B1 is satisfied
as soon as the increment distribution has a bounded $p$th moment. Because, on
the set $\mathsf{R}_x$, $\pi(x+y) \leq \pi(x)$, we similarly have $|\Delta(x)|
\leq \int |y| q(y) \lleb(dy)$ showing, that B2 is
satisfied with $\beta=0$; nevertheless, in some examples, for $\beta=0$,
$\Delta_\infty$ can be zero and the fluid limit model is unstable. In these
cases, it is necessary to use larger $\beta$ (see
Section~\ref{sec:weibulian-tails}).

\subsection{Superexponential target densities}
\label{subsec:Super-exponential target densities}
In this section, we focus on target densities $\pi$ on $\Xset$ which are
\emph{superexponential}. Define $n(x) \eqdef x/|x|$.

\begin{defi}[(\textit{Superexponential p.d.f.})]
A probability density function $\pi$ is said to be \emph{superexponential}
if $\pi$ is positive, has continuous first derivatives and
$\lim_{|x| \to \infty} \pscal{n(x)}{\ell(x)}= - \infty$, where $\ell(x) \eqdef \nabla \log \pi(x)$.
\end{defi}

\begin{figure}[b]

\includegraphics{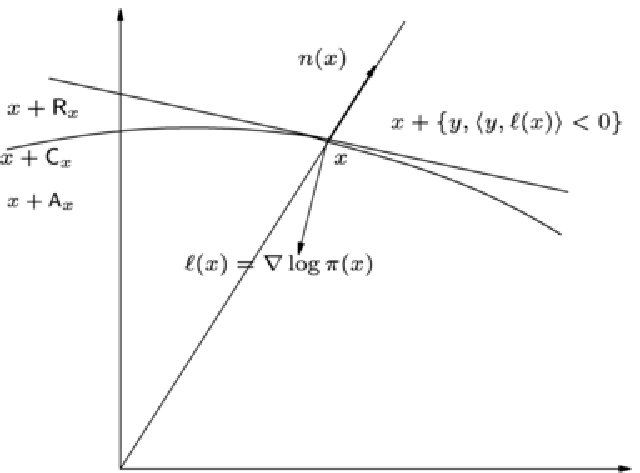}

\caption{}\label{fig:MCMC_principle}
\end{figure}

The condition implies that for any $H > 0$, there exists $R> 0$ such that
\begin{equation}
\label{eq:radial-small}
\frac{\pi(x + an(x))}{\pi(x)} \leq \exp(-a H)\qquad  \mbox{for $|x| \geq R, a \geq  0$},
\end{equation}
that is, $\pi(x)$ is at least exponentially decaying along any ray with the
rate $H$ tending to infinity as $|x|$ goes to infinity. It also implies that
for $x$ large enough, the contour manifold $\mathsf{C}_x \eqdef \{ y \in \Xset,
\pi(x+y) = \pi(x) \}$ can be parameterized by the unit sphere $\unitsphere$
since each ray meets $\mathsf{C}_x$ at exactly one point. In addition, for
sufficiently large $|x|$, the acceptance region $\mathsf{A}_x$ is the set
enclosed by the contour manifold $\mathsf{C}_{x}$ (see Figure~\ref{fig:MCMC_principle}).
Denote by $A \ominus B$ the symmetric difference of the sets $A$ and $B$.

\begin{defi}[(\textit{$q$-radial limit})]
  We say that the family of rejection regions $\{ R_{r x}, r \geq 0, x \in
  \openset\}$ has \emph{$q$-radial limits} over the open cone $\openset \subseteq
  \Xset \setminus \{0\}$ if there exists a collection of sets $\{ R_{\infty,x}, x \in \openset \}$
  such that, for any compact subset $\compactset \subseteq \openset$,
  $\lim_{r \to \infty} \sup_{x \in \compactset} \Q  ( R_{r x} \ominus R_{\infty,x}  ) = 0$.
\end{defi}

\begin{prop}
\label{prop:super-exponential}
Assume that the target density $\pi$ is super-exponential.  Assume, in
addition, that the family $\{ R_{r x}, r \geq 0, x \in \openset \}$ has a $q$-radial
limit over an open cone $\openset \subseteq \Xset \setminus \{0\}$.  Then, for
any compact set $\compactset \subset \openset$,\break $\lim_{r \to \infty} \sup_{x
  \in \compactset} | \Delta(r x)- \Delta_\infty(x)| = 0$, where
$\Delta_\infty(x) \eqdef -\! \int_{R_{\infty,x}} y q(y) \lleb(dy)$.
\end{prop}

The proof is postponed to Section~\ref{proof:sec:superexp}.
The definition of the limiting field $\Delta_\infty$ becomes simple when the rejection region
radially converges to a half-space.

\begin{defi}[(\textit{$q$-regularity in the tails})]
  We say that the target density $\pi$ is \emph{$q$-regular in the tails} over
  $\openset$ if the family $\{ R_{r x}, r \geq 0, x \in \openset \}$ has $q$-radial limits over an open cone
  $\openset \subseteq \Xset \setminus \{ 0 \}$ and there exists a
  continuous function \mbox{$\ell_\infty\dvtx  \Xset \setminus \{0 \} \to \Xset$} such
  that, for all $x \in \openset$,
\begin{equation}
\label{eq:definition-ell-infty}
\Q  \bigl( R_{\infty,x} \ominus  \{y \in \Xset, \pscal{y}{ \ell_\infty(x)} < 0  \}  \bigr) = 0.
\end{equation}
\end{defi}
Regularity in the tails holds with $\ell_\infty(x)= \lim_{r \to \infty} n(\ell(rx))$ when the curvature at $0$ of the contour manifold
$\mathsf{C}_{rx}$ goes to zero as $r \to \infty$; nevertheless, this
condition may still hold in situations where there exists a sequence $\{ x_n \}$
with $\lim |x_n| = \infty$ such that the curvature of the contour
manifolds $C_{x_n}$ at zero can grow to infinity (see Examples
\ref{example:superexp-sharp-wedge}~and~\ref{example:Mixture of Gaussian Densities}).
Assume that
\begin{equation}
\label{eq:definition-q}
q(x) = \mathsf{det}^{-1/2}(\Sigma)   q_0(\Sigma^{-1/2} x),
\end{equation}
where $\Sigma$ is a positive definite matrix and $q_0$ is a \emph{rotationally
  invariant distribution},  that is, $q_0(Ux)= q_0(x)$ for any unitary matrix $U$, and is such that
\[
\int_\Xset y_1^2   q_0(y) \lleb(dy) < \infty.
\]

\begin{prop}
\label{prop:DeltaInftyExpReg}
Assume that the target density $\pi$ is \emph{super-exponential} and
\emph{$q$-regular in the tails} over the open cone $\openset \subseteq \Xset
\setminus \{0\}$.  Then, the SRWM algorithm with proposal $q$ given in
\eqref{eq:definition-q} satisfies assumption \textup{B3} on
$\openset$ with
\begin{equation}
\label{eq:definition-Delta-infty-regular}
\Delta_\infty(x) = m_1(q_0) \frac{\Sigma \ell_\infty(x)}{|\sqrt{\Sigma} \ell_\infty(x) |},
\end{equation}
where $\ell_\infty$ is defined in \eqref{eq:definition-ell-infty} and
$m_1(q_0) \eqdef \int_{\Xset} y_1 \1_{\{ y_1 \geq 0\}} q_0(y) \lleb(dy) >0$, where $y= (y_1,\dots,y_d)$.
\end{prop}

The proof is given in
Section~\ref{proof:sec:superexp}.  If $\Sigma= \mathrm{Id}$ and $\ell_\infty(x)= \lim_{r \to \infty}
n(\ell(rx))$, then the ODE may be seen as a version of steepest ascent algorithm to maximize $\log \pi$. It may appear that
convergence would be faster if $m_1(q_0)$ is increased. While it is true for the ODE, we
cannot reach such a positive conclusion for the algorithm itself because we do
not control the fluctuation of the algorithm around its limit.

\subsubsection{Regular case}

The tail regularity condition and the definition of the ODE limit are more transparent in a class of models which are
very natural in many statistical contexts, namely, the exponential family.
Following \citet{robertstweedie1996}, define the class $\mathcal{P}$ as consisting
of those everywhere positive densities with continuous second derivatives
$\pi$ satisfying
\begin{equation}
\label{eq:regular-exponential}
\pi(x) \propto  g(x) \exp   \{ -p(x)  \},
\end{equation}
where:
\begin{itemize}
\item $g$ is a positive function slowly varying at infinity,  that is, for any $K > 0$,
\begin{equation}
\label{eq:definition-subgeometric}
\limsup_{|x| \to \infty} \inf_{|y| \leq K} \frac{g(x+y)}{g(x)}= \limsup_{|x|
  \to \infty} \sup_{|y| \leq K} \frac{g(x+y)}{g(x)} =1;
\end{equation}

\item $p$ is a positive polynomial in $\Xset$ of even order $m$ and
  $\lim_{|x| \to \infty} p_m(x) = +\infty$, where $p_m$ denotes the polynomial
  consisting only of the $p$'s $m$th order terms.
\end{itemize}

\begin{figure}[b]

\includegraphics{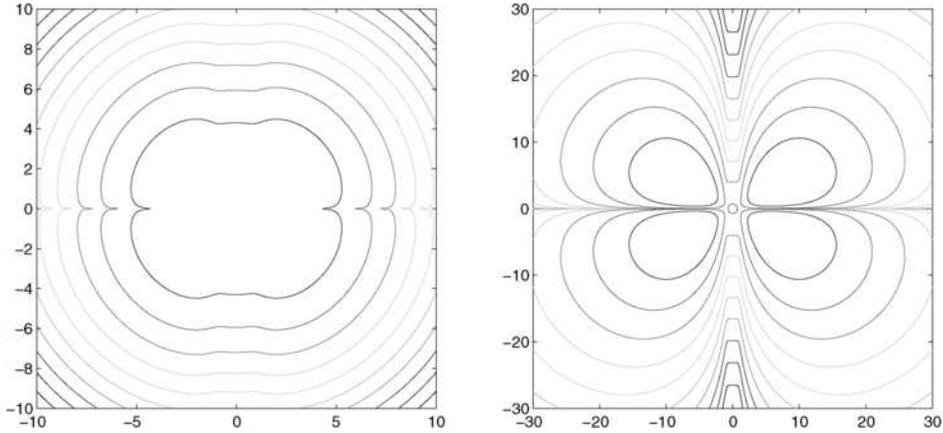}

  \caption{Contour curves of the target densities \protect\eqref{eq:target-super} (left panel)
  and \protect\eqref{eq:target-subexp} with $\delta=0.4$ (right panel).}
  \label{fig:contour_super}
\end{figure}

\begin{prop}\label{prop:stability-fluid-superexponential}
Assume that $\pi \in \mathcal{P}$ and let $q$ be given by
\eqref{eq:definition-q}. Then, $\pi$ is super-exponential, $q$-regular in the tails
over \mbox{$\Xset {\setminus} \{0\}$} with $\ell_{\infty\!}(x)=\break -n  [ \nabla p_{m\!}  ( n(x)  )  ]$.
For any $x \in \Xset \setminus \{0\}$, there exists $T_x > 0$ such that the ODE
$\dot{\mu}= \Delta_\infty(\mu)$ with initial condition $x$ and
$\Delta_\infty$ given by \eqref{eq:definition-Delta-infty-regular} has a unique solution on $[0,T_x)$ and
$\lim_{t \to T_x^-} \mu(t;x)= 0$. In addition, the fluid limit $\QQq{0}{x}$ is deterministic on
$\mathsf{D}([0,T_x],\Xset)$, with support function $\mu(\cdot;x)$.
\end{prop}

The proof is skipped for brevity [see~Fort et al. (\citeyear{FortMeynMoulinesPriouret06})].
Because all the solutions of the initial value problem
$\dot{\mu}= -m_1(q_0) \sqrt{\Sigma}n [\sqrt{\Sigma} \nabla
p_m(n(\mu)) ]$, $\mu(0)= x$ are zero after a fixed amount of time $T$
for any initial condition on the unit sphere, we may apply
Theorem~\ref{theo:stability-fluid-limit-deterministic-smooth}. We have, from
Theorem~\ref{theo:PetiteIrred} and
Theorem~\ref{theo:stability-fluid-limit-deterministic-smooth}, the
following.

\begin{figure}

\includegraphics{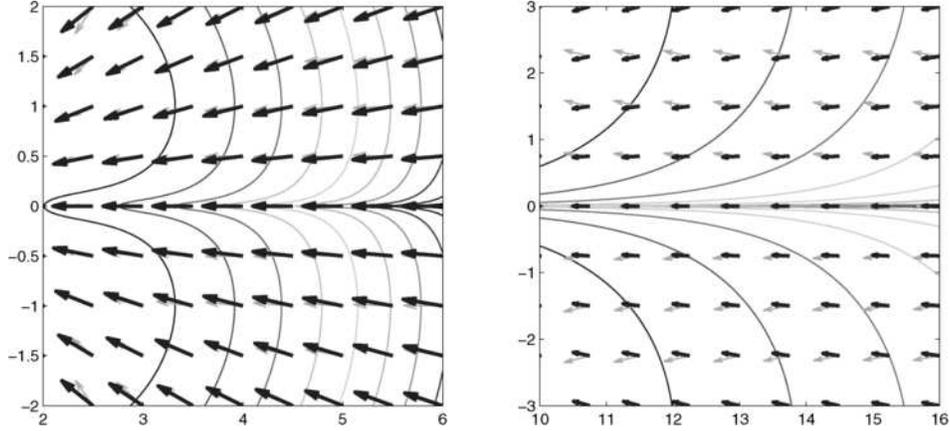}

  \caption{Grey  lines: $\Delta$; Black  lines: $\Delta_\infty$ for the target
  densities~\protect\eqref{eq:target-super} (left panel) and \protect\eqref{eq:target-subexp} with $\delta=0.4$ (right panel). }\label{fig:delta_super}
\end{figure}

\begin{figure}[b]

\includegraphics{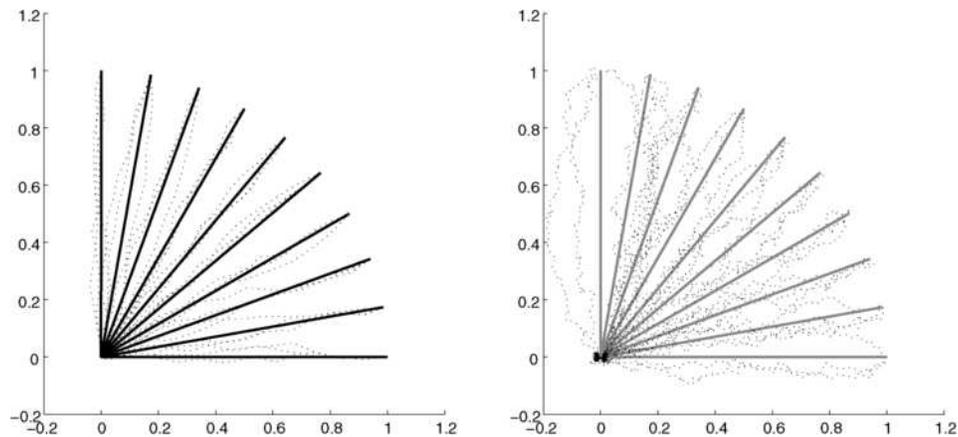}

  \caption{Dotted lines: trajectories of the interpolated process \protect\eqref{eq:InterpolatedProcessIntro} for the
  random walk Metropolis--Hastings (SRWM) algorithm
    for a set of initial conditions on the unit sphere in $(0,\pi/2)$ for
    the target densities~\protect\eqref{eq:target-super} (left panel)
    and~\protect\eqref{eq:target-subexp} (right panel); Solid lines: flow of the
    associated ODE.}
  \label{fig:SuperExp_fluid}
\end{figure}

\begin{theo}
  Consider the SRWM Markov chain with target distribution $\pi \in \mathcal{P}$
  and increment distribution $q$ having a moment of order $p>1$ and
  satisfying~\eqref{eq:definition-q}. Then, for any $1 \leq u \leq p$, the SRWM
  Markov chain is $(f_u,r_u)$-ergodic with
\[
f_u(x) = 1+|x|^{p-u},\qquad  r_u(t) \sim t^{u-1}.
\]
\end{theo}

\begin{figure}

\includegraphics{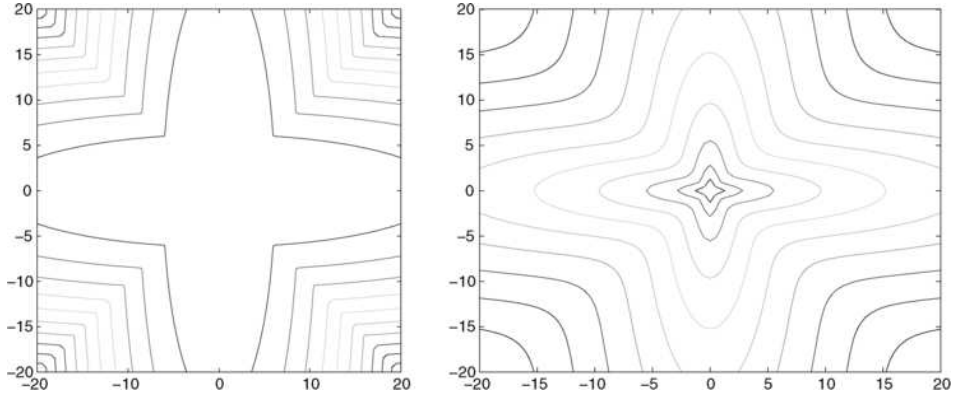}

  \caption{Contour plot of the target densities \protect\eqref{eq:target-mixture-gaussian} (left panel)
  and (\protect\ref{eq:target:mixture-weibul-distribution}) (right panel). } \label{fig:contour_trefle}
\end{figure}

\begin{example}\label{example:superexp-sharp-wedge}
To illustrate our findings, consider the target density, borrowed from
\citet{jarnerhansen2000}, Example~5.3,
\begin{equation}
\label{eq:target-super}
\pi(x_1,x_2) \propto (1 + x_1^2 + x_2^2 + x_1^8 x_2^2) \exp  \bigl( -(x_1^2+x_2^2)  \bigr).
\end{equation}
The contour curves are illustrated in Figure~\ref{fig:contour_super}. They are
almost circular except from some small wedges by the $x$-axis. Due to the
wedges, the curvature of the contour manifold at $(x, 0)$ is $(x^6-1)/x$ and
therefore tends to infinity along the $x$-axis [\citet{jarnerhansen2000}].
Since $\pi \in \mathcal{P}$, Proposition
\ref{prop:stability-fluid-superexponential} shows that $\pi$ is
super-exponential, regular in the tails and $\ell_\infty (x)= -n(x)$. Taking $q
\sim \mathcal{N}(0,\sigma^2 \mathrm{Id})$, $\Delta_\infty(x)= - \sigma n(x) /
\sqrt{2 \pi}$ and the (Carath\'eodory) solution of the initial value problem
$\dot{\mu}= \Delta_\infty(\mu)$,
$\mu(0)= x$ is given by $\mu(t;x)= (|x|- \sigma
t/\sqrt{2\pi}) \1 \{ \sigma t \leq \sqrt{2\pi}|x|\}x/|x|$.  Along the sequence
$\{x_k\eqdef (k,\pm k^{-4}) \}_{k \geq 1}$, the normed gradient $n[\ell(x_k)]$
converges to $(0,\pm 1)$, showing that whereas $\ell_\infty$ is the radial
limit of the normed gradient $n[\ell]$ (i.e.,  for any $u \in \mathsf{S}$,
$\lim_{\lambda \to \infty} n[\ell(\lambda u)]= \ell_\infty(u)$), $\limsup_{|x|
  \to \infty} |n[\ell(x)]-\ell_\infty(x)|=2$. Therefore, the normed gradient
$n[\ell(x)]$ does not have a limit as $|x| \to \infty$ along the $x$-axis.
Nevertheless, the fluid limit exists and is extremely simple to determine.
Hence, the ergodicity of the SRWM sampler with target distribution
(\ref{eq:target-super}) may be established [note that for this example, the
theory developed in~\citet{robertstweedie1996} and
in~\citet{jarnerhansen2000} does not apply].  The functions $\Delta$ and
$\Delta_\infty$ are displayed in Figure~\ref{fig:delta_super}.  The flow of the
initial value problem $\dot{\mu}= \Delta_\infty
(\mu)$ for a set of initial conditions on the unit sphere
in  $(0, \pi/2)$ is displayed in Figure~\ref{fig:SuperExp_fluid}.
\end{example}

\begin{figure}

\includegraphics{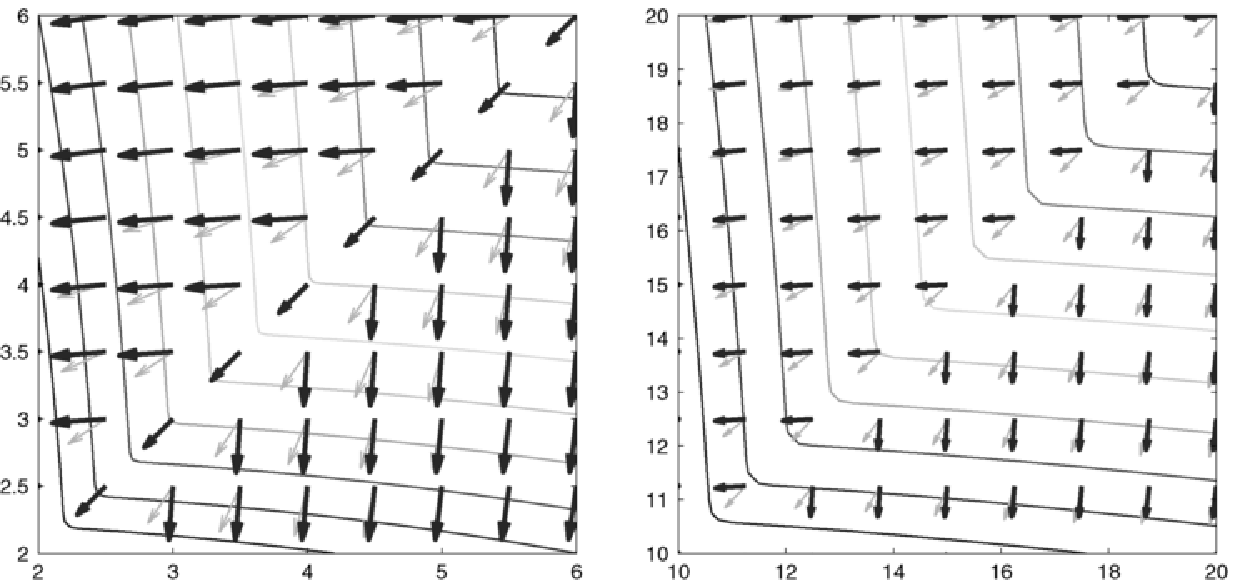}

  \caption{Grey  lines: $\Delta$; Black  lines: $\Delta_\infty$ for the target density \protect\eqref{eq:target-mixture-gaussian}
  (left panel) and
 (\protect\ref{eq:target:mixture-weibul-distribution})    (right panel).} \label{fig:delta_trefle}
\end{figure}

\begin{figure}[b]

\includegraphics{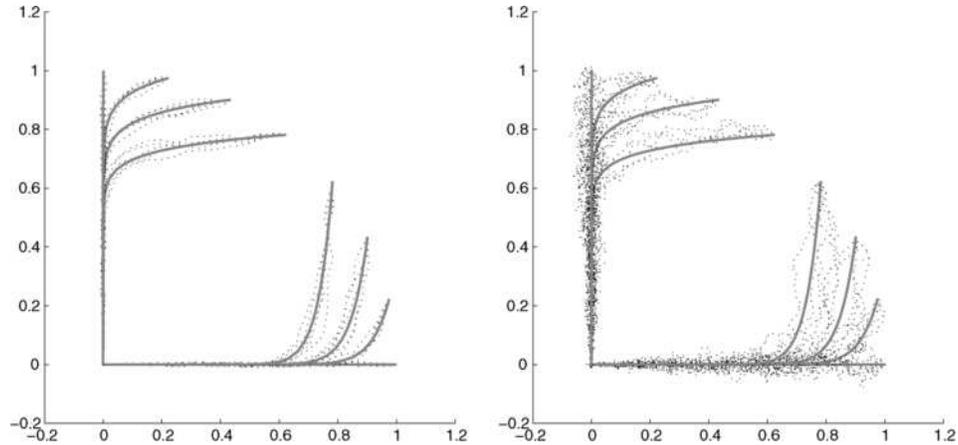}

  \caption{Dotted lines: interpolated process for a set of initial conditions on the
  unit sphere for the target density \protect\eqref{eq:target-mixture-gaussian} (left panel)
    and (\protect\ref{eq:target:mixture-weibul-distribution}) (right panel); Solid lines:
    flow of the initial value problem $\dot{\mu}=
    h(\mu)$ with $h(x) = |x|^{-\beta} \Delta_\infty(x)$; $\beta=
    0$ and $\Delta_\infty$ are given by Lemma \protect\ref{lem:champ-de-trefle} (left
    panel) and $\beta, \Delta_\infty$ are given by
    Lemma~\protect\ref{lem:champ-de-weibuliantrefle} (right panel).}
  \label{fig:trefle-skeleton-1000}
\end{figure}

\subsubsection{Irregular case} We give  an example for which, in
Proposition~\ref{prop:super-exponential}, $\openset \subsetneq \Xset \setminus
\{0\}$.

\begin{example}
\label{example:Mixture of Gaussian Densities}
In this example [also borrowed from \citet{jarnerhansen2000}], we consider the
mixture of two Gaussian distributions on $\rset^2$.  For some $a^2>1$ and $0 <
\alpha < 1$, set
\begin{equation}
\label{eq:target-mixture-gaussian}
\pi(x) \propto    \alpha \exp  (- (1/2) x' \Gamma_1^{-1} x  ) + (1-\alpha)  \exp  ( - (1/2) x' \Gamma_2^{-1} x  ),
\end{equation}
where $\Gamma_1^{-1} \eqdef \mathrm{diag}(a^2,1)$ and $\Gamma_2^{-1} \eqdef
\mathrm{diag}(1,a^2)$.  The contour curves for $\pi$ with $a=4$ are illustrated
in Figure~\ref{fig:contour_trefle}.  We see that the contour curves have some
sharp bends along the diagonals that do not disappear in the limit, even though
the contour curves of the two components of the mixtures are smooth ellipses.
Equation~(51) of \citet{jarnerhansen2000}, indeed shows that the curvature of the
contour curve on the diagonal tends to infinity. As shown in the following
lemma, however, this target density is regular in the tails over $\openset=
\Xset \setminus \{x=(x_1,x_2) \in \rset^2, |x_1| = |x_2|\}$ (and not over
$\Xset \setminus \{0\}$). More precisely, we have the following.

\begin{lem}
\label{lem:champ-de-trefle}
For any $\varepsilon > 0$, there exist $M$ and $K$ such that
\begin{equation}
\label{eq:champ-de-trefle-approx}
\sup_{|x| \geq K,  |   |x_1| - |x_2|    | \geq M}  | \Delta(x) - \Delta_\infty(x)  | \leq \varepsilon,
\end{equation}
where $\Delta_\infty(x) \eqdef -\! \int \1_{R_{\infty, x}}(y)   y q(y) \lleb(dy)$ with
$R_{\infty,x} \eqdef \{y,\pscal{y}{\Gamma_2^{-1} x}   \geq 0 \}$ if $|x_1| > |x_2|$ and
$R_{\infty,x} \eqdef \{y, \pscal{y}{\Gamma_1^{-1} x}   \geq 0 \}$ otherwise.
\end{lem}

The proof is postponed to Section~\ref{sec:proof:MCMC-Ex2}.  Since $q$
satisfies~(\ref{eq:definition-q}), when $\Sigma = \mathrm{Id}$, for any $x \in
\openset$, we have either $\Delta_\infty(x) = - c_q n(\Gamma_2^{-1} x)$ if
$|x_1| > |x_2|$ or $\Delta_\infty(x)= -c_q n(\Gamma_1^{-1} x)$ if $|x_1| <
|x_2|$, where $c_q$ is a constant depending on the increment distribution $q$.
This is illustrated in Figure~\ref{fig:delta_trefle}, which displays the
functions $\Delta$ and $\Delta_\infty$ and shows that these two functions are
asymptotically close outside a band along the main diagonal. The flows of the
initial value problem $\dot{\mu}= \Delta_\infty
(\mu)$ for a set of initial conditions in $(0, \pi/2)$ are
displayed in Figure~\ref{fig:trefle-skeleton-1000}.

\begin{figure}[b]

\includegraphics{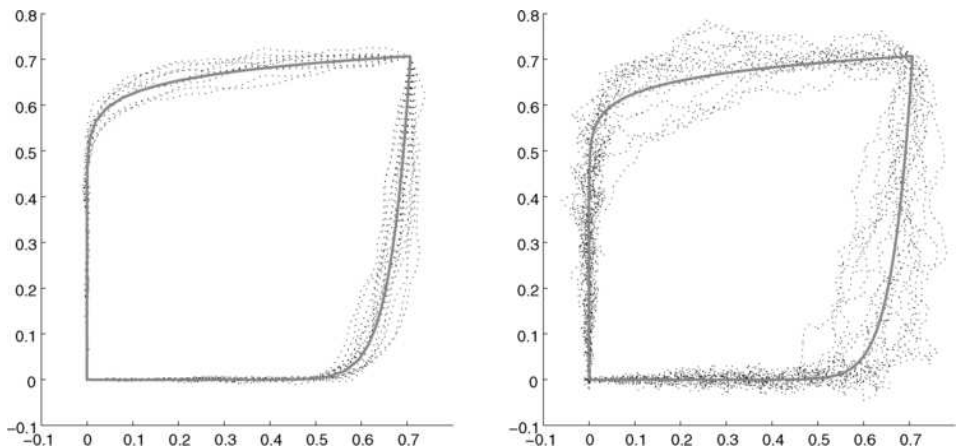}

  \caption{Dotted lines: trajectories of the interpolated process
  \protect\eqref{eq:InterpolatedProcessIntro} for the SRWM with target density
    \protect\eqref{eq:target-mixture-gaussian} (left panel) and
    \protect\eqref{eq:target:mixture-weibul-distribution} (right panel) and initial
    condition $(1/\sqrt{2},1/\sqrt{2})$; Solid lines: flow of the associated
    ODE.}\label{fig:trefle-diago-1000}
\end{figure}

We now prove that
Theorem~\ref{theo:stability-fluid-limit-deterministic-nonsmooth} applies.
Conditions
B1--B2
hold, as discussed above. Condition B3 results from
Lemma~\ref{lem:champ-de-trefle}. It remains to prove that B4
and conditions (i)--(iii)
are verified. The proof of condition (i) is certainly
the most difficult to check in this example.

\begin{prop}
\label{prop:Ex2-cond1PropStab}
Consider the SRWM Markov chain with target distribution given by
\textup{(\ref{eq:target-mixture-gaussian})}. Assume that $q$ is rotationally invariant
and with compact support. Then, \textup{B4} as well as conditions
\textup{(i), (ii)} and \textup{(iii)} of
Theorem~\textup{\ref{theo:stability-fluid-limit-deterministic-nonsmooth}} hold.
\end{prop}

A detailed proof is provided in Section~\ref{sec:proof:MCMC-Ex2}.
Note that the fluid limit model is not deterministic in this example: for $x$
on the diagonal in $\Xset$, the support of the fluid limit $\QQq{0}{x}$
consists of two trajectories, each of which are solutions of the ODE. This is
illustrated in Figure~\ref{fig:trefle-diago-1000}.  By
Theorem~\ref{theo:stability-fluid-limit-deterministic-nonsmooth} and the
discussion above, we may conclude that if the increment distribution $q$ is
compactly supported, then the SRWM Markov chain with target distribution $\pi$ given
by (\ref{eq:target-mixture-gaussian}) is \mbox{$(f_u,r_s)$-}ergodic with $f_u(x) = 1 +
|x|^u$ and $r_s(t) \sim t^s$ for any $u \geq  0$ and $s\geq  0$.
\end{example}

\subsection{Subexponential density}
\label{sec:weibulian-tails}
In this section, we focus on target densities $\pi$ on $\Xset$ which are
\emph{subexponential}. We assume that $q$ satisfies (\ref{eq:definition-q}) and
has moment of order $p \geq 2$. This section is organized as above: we start
with the regular case (Example~\ref{example:subexp-sharp-edge}) and then
consider the irregular case (Example~\ref{example:mixture-weibulian}).

\begin{defi}[({Subexponential p.d.f.})]
  A probability density function $\pi$ is said to be \emph{subexponential} if
  $\pi$ is positive with continuous first derivatives,
  $\pscal{n(x)}{n(\ell(x))} < 0$ for all sufficiently large $x$ and $\lim_{|x|
    \to \infty} | \ell(x) |= 0$.
\end{defi}

The condition implies that for any $R < \infty$, $\lim_{|x| \to \infty}
\sup_{|y| \leq R} \pi(x + y)/\pi(x) = 1$, which implies that $\lim_{|x| \to
  \infty} |\Delta(x)| = 0$.  Subexponential target densities provide examples
that require the use of positive $\beta$ in the normalization to obtain a
nontrivial fluid limit model.

The condition $\pscal{n(x)}{n(\ell(x))} < 0$ for all sufficiently large $|x|$
implies that for $\epsilon$ small enough, the contour manifold
$\mathsf{C}_\epsilon$ can be parameterized by the unit sphere (see the
discussion above) and that for sufficiently large $|x|$, the acceptance region
$\mathsf{A}_x$ is the set enclosed by the contour manifold $\mathsf{C}_{x}$ (see
Figure~\ref{fig:MCMC_principle}).

\begin{defi}[{[$q$-regularity in the tails (subexponential)]}]
  We say that $\pi$ is \emph{$q$-regular in the tails} over an open cone $\openset
  \subseteq \Xset \setminus \{0\}$ if there exists a continuous function
  $\ell_\infty \dvtx  \openset \to \Xset$ and $\beta \in (0,1)$ such that, for any
  compact set $\compactset \subset \openset$ and any $K>0$,
\begin{eqnarray*}
  \lim_{r \to \infty} \sup_{x \in \compactset} \int_{\mathsf{R}_{rx} \cap  \{y, |y| \leq K \} }
  \bigg| r^\beta |x|^\beta   \biggl\{ \frac{\pi(rx+y)}{\pi(rx)}-1   \biggr\}- \pscal{\ell_\infty(x)}{y}  \bigg| q(y)  \lleb(dy)  &=&  0,
  \\
  \lim_{r \to \infty} \sup_{x \in \compactset} \Q  \bigl( \mathsf{R}_{rx}
    \ominus  \{y, \pscal{\ell_\infty(x)}{y} \geq 0  \}  \bigr)&=&  0.
\end{eqnarray*}
\end{defi}

\begin{prop}
\label{prop:subexponential}
Assume that the target density $\pi$ is subexponential and $q$-regular in the tails
over an open cone  $\openset \subseteq \Xset \setminus \{0\}$ and that $q$
satisfies \eqref{eq:definition-q}.  Then, for any compact set $\compactset
\subset \openset$, $\lim_{r \to \infty} \sup_{x \in \compactset} | r^\beta
|x|^\beta \Delta(r x ) - \Delta_\infty( x )| = 0$ with
\[
\Delta_\infty(x) \eqdef \int_{\{y, \pscal{\ell_\infty(x)}{y} \geq 0 \}} y
\pscal{\ell_\infty(x)}{y} q(y) \lleb(dy) = m_2(q_0) \Sigma \ell_\infty(x),
\]
where $m_2(q_0) \eqdef \int_{\Xset} y_1^2 \1 _{\{ y_1 \geq 0 \}} q_0(y)
\lleb(dy) >0$.
\end{prop}

The proof is similar to Proposition \ref{prop:super-exponential} and is omitted
for brevity.  Once again, if the curvature of the contour curve goes to zero at
infinity, then $\ell_\infty(x)$ is, for large $x$, asymptotically colinear to
$n[\nabla \log \pi(x)]$.  However, whereas $|\nabla \log \pi(x)| \to 0$ as $|x|
\to \infty$, the renormalization prevents $\ell_\infty(x)$ from vanishing at
$\infty$; on the contrary, it converges radially to a constant along each ray.
As above, the tail regularity condition may still hold, even when the curvature
goes to infinity; see Example \ref{example:subexp-sharp-edge}. As above, the
subexponential tail regularity condition and the definition of the ODE limit
are more transparent in the Weibullian family.
Mimicking the construction above, define, for $\delta > 0$, the class
$\mathcal{P}_\delta$ as consisting of those everywhere positive densities with
continuous second derivatives $\pi$ satisfying
\begin{equation}
\label{eq:regular-subexponential}
\pi(x) \propto  g(x) \exp   \{ -p^\delta(x)  \},
\end{equation}
where $g$ is a positive function slowly varying at infinity [see
\eqref{eq:definition-subgeometric}] and $p$ is a positive polynomial in $\Xset$ of
even order $m$ with $\lim_{|x| \to \infty} p_m(x) = + \infty$.

\begin{prop}
\label{prop:stability-fluid-subexponential}
Assume that $\pi \in \mathcal{P}_\delta$ for some $0 < \delta < 1/m$ and let
$q$ be given by \eqref{eq:definition-q}. Then, $\pi$ is subexponential and $q$-regular
in the tails with $\beta= 1-m\delta$ and $\ell_\infty(x)= -\delta
p_m^{\delta-1} ( n(x)  ) \nabla p_m ( n(x)  ) $.  For any $x
\in \Xset \setminus \{0\}$, there exists $T_x > 0$ such that the ODE
$\dot{\mu}= h(\mu)$ with initial condition $x$
and $h$ given by
\begin{equation}
\label{eq:definition-h-subexp}
h(x)= -\delta |x|^{-(1-m\delta)} m_2(q_0)  p_m^{\delta-1}(n(x))   \Sigma \nabla p_m(n(x))
\end{equation}
has a unique solution on $[0,T_x)$ and $\lim_{t \to T_x^-} \mu(t;x)=
0$. In addition, the fluid limit $\mathbb{Q}^\beta_x$ is deterministic on
$\mathsf{D}([0,T_x],\Xset)$, with support function
$\mu(\cdot;x)$.
\end{prop}

We may apply Theorem~\ref{theo:stability-fluid-limit-deterministic-smooth}.
From Theorem~\ref{theo:PetiteIrred} and Proposition~\ref{prop:stability-fluid-subexponential} we
have the following.

\begin{theo}
  Consider the SRWM Markov chain with target distribution $\pi$ on
  $\mathcal{P}_\delta$ and increment distribution $q$ having a moment of order
  $p \geq 2$ and satisfying \textup{(\ref{eq:definition-q})}. Then, for any $1 \leq u
  \leq p /(2-m \delta)$, the SRWM Markov chain is \mbox{$(f_u,r_u)$-}ergodic with
\[
f_u(x) = 1+ |x|^{p-u(2-m \delta)},\qquad     r_u(t) \sim t^{u-1}.
\]
\end{theo}

\begin{example}
\label{example:subexp-sharp-edge}
Consider the subexponential Weibullian family derived from
Example~\ref{example:superexp-sharp-wedge},
\begin{equation}
\label{eq:target-subexp}
\pi(x_1,x_2) \propto (1 + x_1^2 + x_2^2 + x_1^8 x_2^2)^{\delta} \exp  \bigl( -(x_1^2+x_2^2)^{\delta}  \bigr).
\end{equation}
The contour curves are displayed in Figure~\ref{fig:contour_super}.  Since $\pi
\in \mathcal{P}_\delta$, Proposition \ref{prop:stability-fluid-subexponential}
shows that $\pi$ is subexponential and regular in the tails with $\beta= 1 - 2
\delta$ and $\ell_\infty(x)= - 2 \delta n(x)$.  Taking $q \sim \mathcal{N}(0,
\sigma^2 \mathrm{Id})$, $\Delta_\infty(x)= - \sigma^2 \delta n(x) $ and the
(Carath\'eodory) solutions of the initial value problem $\dot{\mu}=
|\mu|^{-(1-2\delta)}\Delta_\infty(\mu)$, $\mu(0)= x$
are given by $\mu(t;x)= [|x|^{2(1-\delta)} - 2\sigma^2
\delta(1-\delta)t]^{0.5 (1-\delta)^{-1}} n(x) \times\break \1_{|x|^{2(1-\delta)} - 2\sigma^2
  \delta(1-\delta)t \geq 0}$. Here, again, the gradient $\ell(x)$ (even properly
normalized) does not have a limit as $|x| \to \infty$ along the $x$-axis, but
the fluid limit model is simple to determine.  Hence, the ergodicity of the
SRWM sampler with target distribution (\ref{eq:target-subexp}) may be
established [note that for this example, the theory developed in
\citet{fortmoulines2003SPA} and \citet{doucfortmoulinessoulier2004}
does not apply].  The functions $\Delta$ and $\Delta_\infty$ are displayed in
Figure~\ref{fig:delta_super}.  The flow of the initial value problem
$\dot{\mu}= h(\mu)$ for a set of initial
conditions on the unit sphere in $(0, \pi/2)$ is displayed in
Figure~\ref{fig:SuperExp_fluid}, together with trajectories of the interpolated
process.
\end{example}

\begin{example}
\label{example:mixture-weibulian}
Consider the mixture of bivariate Weibull distributions [see
\citet{patradey1999} for applications],
\begin{eqnarray}\label{eq:target:mixture-weibul-distribution}
\pi(x) &\propto&    \alpha (x' \Gamma_1^{-1} x)^{\delta-1} \exp  (- (1/2) (x' \Gamma_1^{-1} x)^\delta  )
\nonumber
\\[-8pt]
\\[-8pt]
\nonumber
&&\!{} + (1-\alpha) (x' \Gamma_2^{-1} x)^{\delta-1} \exp  ( - (1/2) (x' \Gamma_2^{-1} x)^\delta  ),
\end{eqnarray}
where $\Gamma_i$, $i=\pm 1,2$, are defined in Example~\ref{example:Mixture of
  Gaussian Densities} and $0< \alpha <1$. Similarly to
Example~\ref{example:Mixture of Gaussian Densities}, the curvature of the
contour curve on the diagonal tends to infinity; nevertheless, the target
density is regular in the tails over $\openset= \Xset \setminus \{x=(x_1,x_2)
\in \rset^2, |x_1| = |x_2|\}$. More precisely, we have the following.

\begin{lem}
\label{lem:champ-de-weibuliantrefle}
For any $\varepsilon > 0$, there exist $M$ and $K$ such that
\begin{equation}
\label{eq:champ-de-trefle-approx-weibul}
\sup_{|x| \geq K,  |   |x_1| - |x_2|    | \geq M}  \big| |x|^\beta \Delta(x) - \Delta_\infty(x)  \big| \leq \varepsilon,
\end{equation}
where $\beta \eqdef 1 - 2 \delta$ and $\Delta_\infty(x) \eqdef - m_2(q_0)
|x|^\beta \delta (x'\Gamma_2^{-1}x)^{\delta-1} \Sigma \Gamma_2^{-1} x$ if
$|x_1| > |x_2|$ and $\Delta_\infty(x) \eqdef -m_2(q_0) |x|^\beta \delta
(x'\Gamma_1^{-1}x)^{\delta-1} \Sigma \Gamma_1^{-1} x$ otherwise.
\end{lem}

We can then establish the analog of Proposition~\ref{prop:Ex2-cond1PropStab}
for the target
distribution~(\ref{eq:target:mixture-weibul-distribution}), again
assuming that the proposal distribution $q$ has compact support.  The
details are omitted for brevity. From the discussions above, the SRWM Markov
chain with target distribution $\pi$ given by
(\ref{eq:target:mixture-weibul-distribution}) is $(f_u,r_s)$-ergodic with
$f_u(x) = 1+|x|^u$ and $r_s(t) \sim t^s$ for all $u\geq 0$, $ s \geq 0$.
\end{example}

\section{Conclusions}
ODE techniques provide a general and powerful approach to establishing
stability and ergodic theorems for a Markov chain.  In typical applications, the
assumptions of this paper hold for \textit{any $p>0$} and, consequently, the
ergodic Theorem~\ref{theo:(f,r)-ergodicity} asserts that the mean of
\textit{any} function with polynomial growth converges to its steady-state mean
faster than \textit{any} polynomial rate.  The counterexample presented in
\citet{gamarnikmeyn2005} shows that, in general, it is impossible to obtain a
geometric rate of convergence, even when $\Delta$, $\{\epsilon_k\}$ and the
function $f$ are bounded.

The ODE method developed within the queueing networks research community has
undergone many refinements and has been applied in many very different
contexts.  Some of these extensions might serve well in other applications,
such as MCMC.  In particular, we should note the following points.

\begin{longlist}[(ii)]
\item[(i)] Control variates have been proposed previously in MCMC to speed
  convergence and construct stopping rules [\citet{robert1998}].  The fluid
  model is a convenient tool for constructing control variates for application
  in the simulation of networks.  The resulting simulators show dramatic
  performance improvements in numerical experiments: a hundredfold variance
  reduction is obtained in experiments presented in \citet{hendersonmeyn1997}
  and \citet{hendersonmeyntadic2003} based on marginal additional
  computational effort.  Moreover, analytical results demonstrate that the
  asymptotic behavior of the controlled estimators are greatly improved
  [Meyn (\citeyear{meyn2005,meyn2006,meyn-book2006})].  It is likely that both the theory
  and methodology can be extended to other applications.

\item[(ii)] A current focus of interest in the networks community is the reflected
  diffusion model obtained under a ``heavy traffic scaling.''  An analog of
  ``heavy-traffic'' in MCMC is the case $\beta>0$ considered in this paper; the
  larger scaling is necessary to obtain a nonstatic fluid limit (see Theorem
  \ref{theo:weak-compacity}).  We have maintained $\beta<1$ in order  to obtain a
  deterministic limit.  With $\beta = 1$, we expect that a diffusion limit will
  be obtained for the scaled MH algorithm under general conditions. This will
  be an important tool in the subexponential case.  In the fluid setting of
  this paper, when $\beta>0$, it is necessary to assume a great deal of
  regularity on the densities $\pi$ and $q$ appearing in the MH algorithm to
  obtain a meaningful fluid limit model.  We expect that very different
  regularity assumptions will be required to obtain a diffusion limit and that
  new insights will be obtained from properties of the resulting diffusion
  model.
\end{longlist}

\section{Proofs of the main results}
\label{sec:Proofs}

\subsection{State-dependent drift conditions}
\label{sec:proof:theo:(f,r)-ergodicity}
In this section, we improve the state-dependent drift conditions proposed by
\citet{filonov1989} for discrete state space and later extended by
\citet{meyntweedie1994d} for general state space Markov chains [see also
\citet{meyntweedie1993} and \citet{robert2000} for additional references and
comments].

Following \citet{nummelintuominen1983}, we denote by $\Lambda$ the
set of nondecreasing sequences $r = \{r(n)\}_{n \in \nset}$
satisfying $\lim_{n
  \to \infty} \downarrow \log r(n)/n = 0$,  that is, $\log r(n)/n$ converges
to zero monotonically from above.  A sequence $r \in \Lambda$ is said
to be \textit{subgeometric}.  Examples include polynomial sequences $r(n) =
(n+1)^\delta$ with $\delta > 0$ and \emph{truly} subexponential sequences,
$r(n) = (n+1)^\delta e^{c n^\gamma}$ [$c > 0$ and $\gamma \in (0,1)$].
Denote by $\mathcal{C}$ the set of functions
\begin{eqnarray}\label{eq:definitionsetC}
&& \mathcal{C} \eqdef  \biggl\{ \phi\dvtx  [1,\infty) \to \rset^+,
\mbox{$\phi$ is concave, monotone nondecreasing,}
\nonumber
\\[-8pt]
\\[-8pt]
\nonumber
&&\hspace*{18.5mm}  \mbox{differentiable and }      \inf_{\{ v \in [1,\infty)\}} \phi(v) > 0, 
   \lim_{v\to\infty}\phi'(v) = 0  \biggr\}.
\end{eqnarray}
For $\phi \in \mathcal{C}$, define $\Hphi(v) \eqdef \int_{1}^v (1/\varphi(x))\,dx$.
 The function $\Hphi \dvtx  [1,\infty) \to [0,\infty)$ is increasing and
$\lim_{v \to \infty} \Hphi(v) = \infty$; see [\citet{doucfortmoulinessoulier2004},
Section~2].  Define, for $u \geq 0$, $ r_\varphi(u)
\eqdef \varphi \circ \Hphi^{-1}(u) / \varphi \circ \Hphi^{-1}(0)$, where
$\Hphi^{-1}$ is the inverse of $\Hphi$. The function $u \mapsto r_\varphi(u)$
is log-concave and thus the sequence $\{ r_\varphi(k) \}_{k \geq 0}$ is
subgeometric.  Polynomial functions $\varphi(v) = v^\alpha$, $\alpha \in
(0,1)$ are associated with polynomial sequences $r_\varphi(k)
=(1+(1-\alpha)k)^{\alpha/(1-\alpha)}$.

\begin{prop}
\label{prop:StateDep}
Let $f\dvtx \Xset \to [1, \infty)$ and $ V\dvtx  \Xset \to [1, \infty)$ be measurable
functions, $\varepsilon \in (0,1)$ be a constant and $C \in \Xsigma$ be a set.
Assume that $\sup_C f/V <\infty$ and that there exists a stopping time
$\tau\geq 1$ such that, for any $x \notin C$,
\begin{equation}
\PE_x \Biggl[\sum_{k=0}^{\tau-1} f(\Phi_k) \Biggr] \leq
V(x)\quad\mbox{and}\quad
\PE_x [V(\Phi_{\tau})  ] \leq (1-\varepsilon) V(x).  \label{eq:CS_ContMomTauC-gen}
\end{equation}
Then, for all $x \notin C$, $\PE_x [ \sum_{k=0}^{\tau_C } f(\Phi_k)
   ] \leq  ( \varepsilon^{-1} \vee \sup_C f/V  )  V(x)$. If,\break in addition, we assume that $\sup_{x \in C} \{ f(x) +
  \PE_x[V(\Phi_1)]\}<\infty$, then\break $\sup_{x \in C} \PE_x [ \sum_{k=0}^{\tau_C } f(\Phi_k)  ] < \infty$.
\end{prop}

\begin{pf}
  Set $\tilde \tau \eqdef \tau \1_{C^c}(\Phi_0) + \1_C(\Phi_0)$ and define
  recursively the sequence $\{ \tau^{n} \}$ by $\tau^{0} \eqdef 0$, $\tau^{1}
  \eqdef \tilde \tau$ and $\tau^{n} \eqdef \tau^{n-1} + \tilde \tau \circ
  \theta^{\tau^{n-1}}$, where $\theta$ is the shift operator.  For any $n \in
  \nset$, define by $\bphi_n= \Phi_{\tau^{n}}$ the chain sampled at the
  instants $\{ \tau^{n} \}_{n \geq 0}$.  $\{ \bphi_n \}_{n \geq 0}$ is a Markov
  chain with transition kernel $\bP(x,A) \eqdef \PP_x(\Phi_{\tilde \tau} \in
  A)$, $x \in \Xset$, $A \in \Xsigma$. Equation \eqref{eq:CS_ContMomTauC-gen}
  implies that
\begin{equation}
    \label{eq:EquInterm}
    \bP V(x) = \PE_x [ V(\Phi_{\tilde \tau})  ] \leq V(x) -   F(x)\qquad\mbox{for all } x \notin
    C,
\end{equation}
where $F(x) \eqdef \varepsilon \PE_x  [ \sum_{k=0}^{\tilde \tau-1}  f(\Phi_k)  ]$.
Let $\bar{\tau}_C \eqdef \inf \{n \geq 1, \bphi_n \in C \}$.  Applying the Markov property and the bound $\tau_C \leq \tau^{\bar \tau_C}$,
we obtain, for all $x \notin C$,
\begin{eqnarray*}
   \PE_x \Biggl[ \sum_{k=0}^{\tau_C} f(\Phi_k)   \Biggr] 
   &\leq & \PE_x  \Biggl[ \sum_{k=0}^{\bar \tau_C-1}  \sum_{j=0}^{\tilde \tau  \circ \theta^{\tau^k}-1}f(\Phi_{j+\tau^k})   \Biggr]
     + \PE_x  \bigl[f(\Phi_{\tau^{\bar \tau_C}})   \1_{\{ \tau^{\bar \tau_C } < \infty\} } \bigr]
     \\
   & \leq &\varepsilon^{-1} \PE_x  \Biggl[ \sum_{k=0}^{\bar \tau_C-1} F(\bphi_k)
    \Biggr] +  \biggl( \sup_C \frac{f}{V}  \biggr) \PE_x
    \bigl[V(\Phi_{\tau^{\bar \tau_C}}) \1_{\{ \tau^{\bar \tau_C } < \infty\} }
    \bigr].
\end{eqnarray*}
Furthermore, (\ref{eq:EquInterm}) and the comparison theorem [\citet{meyntweedie1993}, Theorem~11.3.2] applied to the sampled
chain $\{\bphi_n \}_{n \geq 0}$ yields
\[
 \PE_x  \Biggl[ \sum_{k=0}^{\bar \tau_C-1}  F(\bphi_k)   \Biggr]
 + \PE_x  \bigl[V(\Phi_{\tau^{\bar \tau_C}}) \1_{\{\tau^{\bar \tau_C} < \infty \}}
 \bigr]   \leq V(x), \qquad    x \notin C,
\]
which concludes the proof of the first claim. The second claim follows by
writing, for $x \in C$,
\begin{eqnarray*}
  \PE_x \Biggl[ \sum_{k=0}^{\tau_C} f(\Phi_k)  \Biggr]
  & \leq &  2 \sup_C f + \PE_x \Biggl[ \1 { \{X_1 \notin C\}}    \sum_{k=1}^{\tau_C} f(\Phi_k)   \Biggr]
  \\
  &\leq & 2 \sup_C f + \PE_x \Biggl[ \1 { \{X_1 \notin C\}}   \PE_{X_1} \Biggl[  \sum_{k=0}^{\tau_C-1} f(\Phi_k)   \Biggr] \Biggr]
  \\
  &\leq&  2 \sup_C f  +  \biggl( \varepsilon^{-1} \vee \sup_C f/V  \biggr)   \PE_x [ \1 {\{ X_1 \notin C \}}    V(X_1)  ].
\end{eqnarray*}\upqed
\end{pf}

\begin{prop}
\label{prop:RateDep} Assume that the conditions of
Proposition~\textup{\ref{prop:StateDep}}\break are satisfied with $f (x)= \phi \circ V(x)$ for
$x \notin C$ with $\phi \in \mathcal{C}$.  Then, for $x \notin C$, $\PE_x
 [\sum_{k=0}^{\tau_C -1} r_{\tilde{\phi}}(k)  ] \leq M^{-1} V(x)$ and
$\sup_{x \in C} \PE_x  [\sum_{k=0}^{\tau_C -1} r_{\tilde{\phi}}(k)  ]<
\infty$, where,
\begin{equation}
\label{eq:definition-tilde-phi}
\mbox{for all $t$\qquad   $\tilde{\phi}(t)\eqdef \phi(Mt)$}\quad \mbox{and}\quad
M  \eqdef \biggl[\varepsilon^{-1} \vee \sup_C \phi \circ V / V\biggr]^{-1}.
\end{equation}
\end{prop}

\begin{pf}
  It is known that $U(x) \eqdef \PE_x  [ \sum_{k=0}^{\sigma_C} \phi \circ
    V(\Phi_k) ]$, where $\sigma_C \eqdef \inf \{ k \geq 0, \Phi_k \in C
  \}$, solves the equations $PU(x) = U(x) - \phi \circ V( x)$, $x \notin C$ and
  $U(x) = \phi \circ V(x)$, $x \in C$ [see~\citet{meyntweedie1993},
  Theorem~14.2.3].  By Proposition~\ref{prop:StateDep}, $U(x) \leq
  M^{-1}   V(x)$ for all $x \notin C$. Hence,
\begin{equation}
  \label{eq:DriftSousGeom}
PU(x) \leq U(x) -    \tilde \phi \circ U(x),\qquad   x \notin C.
\end{equation}
From \eqref{eq:DriftSousGeom} and \citet{doucfortmoulinessoulier2004}, Proposition~2.2,
$\PE_x  [\sum_{k=0}^{\tau_C -1} r_{\tilde \phi}( k)  ] \leq U(x) \leq
M^{-1}  V(x)$, for $x \notin C$. The proof is concluded by noting that for $x \in C$,
\begin{eqnarray*}
\PE_x  \Biggl[\sum_{k=0}^{\tau_C -1}  r_{\tilde \phi}( k)   \Biggr]
&\leq& r_{\tilde \phi}(0) + \PE_x  \Biggl[\1 { \{ \Phi_1 \notin C \}}   \sum_{k=1}^{\tau_C -1}  r_{\tilde  \phi}( k)   \Biggr]
\\
& \leq&  r_{ \tilde \phi}(0) + M^{-1}   \sup_{x \in C} PV(x)< \infty.
\end{eqnarray*}\upqed
\end{pf}

\begin{theo}
\label{theo:state-dependent-drift}
Suppose that $\{ \Phi_n \}_{n \geq 0}$ is a phi-irreducible and aperiodic Markov chain.
Assume that there exist a function $\phi \in \mathcal{C}$,  a measurable function $V\dvtx  \Xset \to [1,\infty)$,
a stopping time $\tau \geq 1$, a constant $\varepsilon \in (0,1)$ and a petite set $C \subset \Xsigma$
such that
  \begin{eqnarray}
    \PE_x \Biggl[\sum_{k=0}^{\tau-1} \phi \circ V(\Phi_k) \Biggr] &\leq&  V(x),\hspace*{11.28mm}  \qquad  x \notin C\label{eq:CS_ContMomTauC1},
    \\
    \PE_x [V(\Phi_{\tau})  ] &\leq& (1-\varepsilon) V(x),\qquad  x \notin C,\label{eq:CS_ContMomTauC2}
    \\
    \sup_C \{ V + P V \} &<& \infty.       \label{eq:CS_ContMomTauC3}
  \end{eqnarray}
$P$ is then positive Harris recurrent with invariant probability $\pi$
and:
\begin{longlist}[(3)]
\item[(1)] \label{theo:state-dependent-drift:item-1} for all $x \in \Xset$, $\lim_{n
    \to \infty} r_{\tilde{\phi}}(n) \tvnorm{P^n(x,\cdot) - \pi} = 0$, where $\tilde{\phi}$ is defined in~\eqref{eq:definition-tilde-phi};

\item[(2)] for all $x \in \Xset$, $\lim_{n \to \infty} \fnorm{P^n(x,\cdot) -
    \pi}{\phi \circ V} = 0$;

\item[(3)] \label{theo:state-dependent-drift:item-2} the fundamental kernel $Z$ is a bounded linear
transformation from $L_\infty^{\phi \circ V}$ to~$L_\infty^V$.
\end{longlist}
\end{theo}

\begin{pf}
  (1--2) By~\citeauthor{tuominentweedie1994} [(\citeyear{tuominentweedie1994}), Theorem~2.1], it is sufficient
  to prove that
\[
\sup_{x \in C} \PE_x  \Biggl[ \sum_{k=0}^{\tau_C -1} r_{\tilde{\phi}}(k)  \Biggr] <
\infty,  \qquad \sup_{x \in C} \PE_x  \Biggl[ \sum_{k=0}^{\tau_C -1} \phi \circ
  V(\Phi_k)  \Biggr] < \infty
\]
and, for all $x \in \Xset$,
\[
\PE_x  \Biggl[ \sum_{k=0}^{\tau_C -1} r_{\tilde{\phi}}(k)  \Biggr] <
\infty,\qquad
\PE_x  \Biggl[ \sum_{k=0}^{\tau_C -1} \phi \circ V(\Phi_k)  \Biggr] < \infty.
\]
In Proposition~\ref{prop:RateDep} we show that the stated assumptions imply such bounds.

(3) By~\citet{glynnmeyn1996}, Theorem~2.3, it is sufficient to prove
that there exist constants $b,c< \infty$ such that for all $x \in \Xset$,
$PW(x) \leq W(x) - \phi \circ V(x) + b \1_C(x)$, with  $W(x) \leq c V(x)$.
This follows from Proposition~\ref{prop:StateDep}, which shows that
$\sup_{x \in C} \PE_x  [ \sum_{k=0}^{\tau_C } \phi \circ V(\Phi_k)  ] < \infty$
and $\PE_x  [ \sum_{k=0}^{\tau_C } \phi \circ V(\Phi_k)  ] \leq c V(x)$ for all $x \notin C$.
\end{pf}

Using an interpolation technique, we derive a rate of convergence  associated with
some $g$-norm, $0 \leq g \leq \phi \circ V$.

\begin{coro}[({Theorem} \ref{theo:state-dependent-drift})]
\label{coro:state-dependent-drift}
For any pair of functions $(\alpha, \beta)$ satisfying $\alpha(u) \beta(v) \leq u+  v$, for all $(u,v) \in \rset^+ \times \rset^+$ and all $x \in
\Xset$,
\[
\lim_n \ \alpha (r_{\tilde{\phi}}(n) )  \fnorm{P^n(x,\cdot) - \pi }{\beta (\phi \circ
  V) \vee 1}=0.
\]
\end{coro}

A pair of functions $(\alpha,\beta)$ satisfying this condition can be constructed by using  Young's
inequality [\citet{krasnoselskiirutickii1961}].

\subsection[Proof of Theorem 1.2]{Proof of Theorem \textup{\protect\ref{theo:weak-compacity}}}\label{sec:proof:theo:weak-compacity}

We preface the proof with a preparatory lemma. For any process $\{ \epsilon_k
\}_{k \geq 1}$, define
\begin{equation}
\label{eq:definition-M-infty}
M_\infty(\epsilon,n) \eqdef \sup_{ 1 \leq l \leq n  } \Bigg| \sum_{k=1}^{l} \epsilon_k  \Bigg|.
\end{equation}

\begin{lem}
\label{lem:ControleFluctuationPhi}
Assume \textup{B1} and \textup{B2}.
\begin{longlist}[(iii)]
\item[(i)]  For all $\kappa>0$, $J$ and $K$ integers with $J < K$,
\begin{eqnarray*}
&& \sup_{0\leq k \leq k+j \leq K, 0 \leq j \leq  J}  | \Phi_{k+j} - \Phi_k  |
\\
&&\qquad
\leq 8 M_\infty(\epsilon,K) + 2 N(\beta,\Delta) \kappa^{-\beta} J + N(\beta,\Delta) + 2 \kappa,
\end{eqnarray*}
where $N(\beta,\Delta)$ is given in \textup{B2}.

\item[(ii)] \label{item2:ControleFluctuationPhi} For all $0 \leq \alpha \leq \beta$
  and all $T>0$, there exists $M$ such that
\[
\lim_{r \to \infty}   \sup_{x \in \Xset} \PP_x  \biggl( \sup_{0 \leq k \leq k+j
    \leq \lfloor T r^{1+\alpha} \rfloor} | \Phi_{k+j} - \Phi_k| \geq M r
 \biggr)=0.
\]

\item[(iii)] \label{item3:ControleFluctuationPhi} For all $T > 0$ and $\varepsilon >
  0$, there exists $\delta > 0$ such that
\[
\lim_{r \to \infty}   \sup_{x \in \Xset} \PP_x  \biggl( \sup_{0 \leq k \leq k+j
    \leq \lfloor T r^{1+\beta} \rfloor, 0 \leq j \leq \lfloor \delta
    r^{1+\beta} \rfloor} | \Phi_{k+j} - \Phi_k| \geq \varepsilon r  \biggr)=0.
\]
\end{longlist}
\end{lem}

\begin{pf} (i) Let $0 \leq j \leq J$ and $0 \leq k \leq K-j$. On the set $ \bigcap_{l=0}^{j-1} \{|\Phi_{k+l}| > \kappa \}$,
\begin{eqnarray}
\label{eq:borne-1}
\nonumber | \Phi_{k+j} - \Phi_k   |  &=&    \Bigg|\sum_{l=k}^{k+j-1} \{ \Phi_{l+1} - \Phi_{l}  \}  \Bigg|
\\
&\leq&   \Bigg| \sum_{l=k+1}^{k+j} \epsilon_{l} \Bigg| + \sum_{l=k}^{k+j-1}  | \Delta(\Phi_{l})  |
\nonumber
\\[-8pt]
\\[-8pt]
\nonumber
&\leq & \Bigg| \sum_{l=k+1}^{k+j} \epsilon_{l} \Bigg| + \sum_{l=k}^{k+j-1} |\Phi_{l}|^{-\beta} N(\beta, \Delta)
\\
\nonumber &\leq& 2 M_\infty(\epsilon, K)
+ J \kappa^{-\beta} N(\beta, \Delta).
\end{eqnarray}
Consider now the case when $| \Phi_{k+l} | \leq \kappa $ for some $0 \leq l \leq j-1$. Define
\[
\tau_{j} \eqdef \inf \{   0 \leq l \leq j-1, |\Phi_{k+l}| \leq \kappa  \}
\]
and
\[
\sigma_{j} \eqdef \sup \{ 0 \leq l \leq j-1, |\Phi_{k+l}| \leq \kappa \} +
1,
\]
which are, respectively, the first hitting time and the last exit time before $j$
of the ball of radius $ \kappa  $. Write
$\Phi_{k+j} - \Phi_k = (\Phi_{k+j} - \Phi_{k+\sigma_{j}}) + (\Phi_{k+\sigma_{j}} - \Phi_{k+\tau_{j}}) + (\Phi_{k+\tau_{j}} - \Phi_{k})$
and consider the three terms separately.  The first term  is
nonnull if $\sigma_{j} <j$; hence,
\begin{eqnarray*}\label{eq:borne-2}
|\Phi_{k+j} - \Phi_{k+\sigma_{j}} | &\leq&   \Bigg|\sum_{l=k+\sigma_{j}+1}^{k+j} \epsilon_{l}  \Bigg| + \sum_{l=k+\sigma_{j}}^{k+j-1}  | \Delta(\Phi_{l}) |
\\
&\leq& 2 M_\infty(\epsilon, K) + J \kappa^{-\beta} N(\beta,\Delta)
\end{eqnarray*}
since, by the definition of $\sigma_{j}$, $|\Phi_{k+l}| > \kappa $ for all
$\sigma_{j} \leq l \leq j-1$. Similarly, for the third term,
\begin{eqnarray}  \label{eq:borne-3}
  | \Phi_{k+\tau_{j}} - \Phi_{k} | &\leq&  \Bigg| \sum_{l=k+1}^{k+\tau_{j}} \epsilon_{l}  \Bigg| + \sum_{l=k}^{k+\tau_{j}-1}  | \Delta(\Phi_{l})  |
  \nonumber
  \\[-8pt]
  \\[-8pt]
  \nonumber
  &\leq& 2 M_\infty(\epsilon,K) + J \kappa^{-\beta} N(\beta,\Delta)
\end{eqnarray}
since, by the definition of $\tau_{j}$, $|\Phi_{l}| > \kappa $ for all $ 0 \leq l < \tau_{j}$. Finally, the second term is bounded by
\begin{eqnarray*}\label{eq:borne-4}
|\Phi_{k+\sigma_{j}} - \Phi_{k+\tau_{j}}| &\leq& |\Phi_{k+\sigma_{j}} - \Phi_{k+\sigma_{j}-1}| + |\Phi_{k+\sigma_{j}-1}| + |\Phi_{k+\tau_{j}}|
\\
&\leq& N(\beta,\Delta) + 2 M_\infty(\epsilon,K) + 2 \kappa.
\end{eqnarray*}
Combining the inequalities above yields the desired result.

(ii) From the previous inequality applied with
$\kappa= \ell r>0$ and $K=J= \lfloor T r^{1+\alpha} \rfloor$, it holds
that
\begin{eqnarray*}
&&  \PP_x  \biggl( \sup_{0 \leq k \leq k+j \leq \lfloor T r^{1+\alpha} \rfloor} |
    \Phi_{k+j} - \Phi_k| \geq 4 M r \biggr)
\\
&&\qquad    \leq 4^p M^{-p} r^{-p}   \sup_{x
    \in
    \Xset} \PE_x  [M_\infty^p(\epsilon, \lfloor T r^{1+\alpha} \rfloor)  ]
    \\
&&\qquad\quad {}  + \1\{N(\beta,\Delta) \geq M r\} + \1 \{2 N(\beta,\Delta) T \geq \ell^{\beta}
  M r^{-\alpha +\beta} \} + \1 \{2 \ell \geq M \}.
\end{eqnarray*}
By Lemma~\ref{lem:borne-M-infty}, the expectation tends to zero uniformly for
$x \in \Xset$. The second term tends to zero when $r \to \infty$. The remaining
two terms are zero with $\ell$ and $M$ chosen so that $\ell^{1+\beta} >
N(\beta,\Delta) T$ and $M>2\ell$.

(iii) The proof follows similarly upon setting $K= \lfloor T r^{1+\beta} \rfloor$, $J= \lfloor \delta r^{1+\beta} \rfloor$ and $\kappa= \ell r$.
\end{pf}

\begin{pf*}{Proof of Theorem~\ref{theo:weak-compacity}}
  Let $\alpha \leq \beta$. A sequence of probability measures on
  $\mathsf{D}(\rset^+,\Xset)$ is said to be
  $\mathsf{D}(\rset^+,\Xset)$-\emph{tight} if it is tight in
  $\mathsf{D}(\rset^+,\Xset)$ and if every weak limit of a subsequence is
  continuous.  By [\citet{billingsley1999}, Theorem 13.2, (13.7), page~140 and
  Corollary,  page~142], the sequence of probability measures $\{
  \mathbb{Q}_{r_n;x_n}^{\alpha} \}_{n \geq 0}$ is
  $\mathsf{C}(\rset_+,\Xset)$-tight if (a) $\lim_{a \to \infty}
  \limsup_n \mathbb{Q}_{r_n;x_n}^{\alpha} \{ \eta \dvtx  |\eta(0)|
  \geq a \}= 0$, (b) $\limsup_{n \to \infty} \mathbb{Q}_{r_n;x_n}^{\alpha} \{
  \eta \dvtx  \sup_{0 \leq t \leq T} |\eta(t)
  -\eta(t-)| \geq a \}=0$ and (c) for all $\kappa>0$ and
  $\varepsilon > 0$, there exist $\delta \in (0,1)$ such that $\limsup_n
  \mathbb{Q}_{r_n;x_n}^{\alpha}  \{ \eta \dvtx  w(\eta,\delta)
    \geq \varepsilon  \} \leq \kappa$, where $w(\eta,\delta)
  \eqdef \sup_{0 \leq s \leq t \leq T, |t-s| \leq \delta} |\eta(t)
  - \eta(s)|$.  Properties (a)--(c) follow immediately from Lemma
  \ref{lem:ControleFluctuationPhi}.  Choose $\alpha < \beta$. Let $\{r_n \}$
  and $\{x_n\}$ be sequences such that $\lim_n r_n = \infty$ and $\lim_n x_n =
  x$.  Let $\varepsilon>0$. We have, for all $n$ sufficiently large that $|x_n -
  x| \leq \varepsilon/2$,
\[
  \PP_{r_n x_n}  \biggl( \sup_{0 \leq t \leq T} |\eta_{r_n}^{\alpha}({t};{x_{n}}) - x | \geq \varepsilon  \biggr)
  \leq \PP_{r_n x_n}  \biggl( \sup_{0 \leq k \leq \lfloor T r_n^{1+\alpha} \rfloor } |\Phi_{k} - r_n x_n| \geq (\varepsilon/2) r_n   \biggr)
  \]
  and we have (b), again by
  Lemma~\ref{lem:ControleFluctuationPhi}(ii).
\end{pf*}

\subsection[Proof of Theorem 1.4]{Proof of Theorem \textup{\protect\ref{theo:(f,r)-ergodicity}}}
\label{sec:proof_theo-frErgo}
We preface the proof by establishing a uniform integrability condition for the
martingale increment sequence $\{ \epsilon_k \}_{k \geq 1}$ and then for the
Markov chain $\{ \Phi_k \}_{k \geq 0}$.

\begin{lem}
\label{lem:uniforme-integrability-martingale}
Assume \textup{B1}. Then, for all $T > 0$,
\begin{eqnarray}
\label{eq:uniforme-integrability-martingale}
\hspace*{10mm} && \lim_{b \to \infty} \sup_{|x|\geq 1} |x|^{-p} \PE_x[ M_\infty^p(\epsilon,\lfloor T |\Phi_0|^{1+\beta} \rfloor) \1 \{ M_\infty(\epsilon,\lfloor T |\Phi_0|^{1+\beta} \rfloor) \geq b |\Phi_0| \} ]
\nonumber
\\[-16pt]
\\
\nonumber
&&\qquad = 0.
\end{eqnarray}
\end{lem}

\begin{pf}
Set $T_{\Phi_0} \eqdef \lfloor T |\Phi_0|^{1+\beta} \rfloor$. For $K \geq 0$, set $\bepsilon{K}{k} \eqdef \epsilon_k \1 \{ |\epsilon_k| \leq K \}$ and $\tepsilon{K}{k} \eqdef \epsilon_k \1 \{ |\epsilon_k| \geq K \}$.
By Lemma
\ref{lem:elem-1}, there exists a constant $C$ (depending only on $p$) such that
\begin{eqnarray*}
&& \PE_x[ M_\infty^p(\epsilon,T_{\Phi_0}) \1 \{ M_\infty(\epsilon,T_{\Phi_0}) \geq b |\Phi_0| \} ]
\\
&&\qquad \leq C \PE_x [ M_\infty^p(\bepsilon{K}{},T_{\Phi_0}) \1 \{ M_\infty(\bepsilon{K}{},T_{\Phi_0}) \geq (b/2) |\Phi_0| \}]  + C \PE_x[M_\infty^p(\teps{K},T_{\Phi_0})].
\end{eqnarray*}
Consider the first term on the right-hand side of the previous inequality.
Using Lemma \ref{lem:elem-2} with $a > 1 \vee 2/p$ and Lemma
\ref{lem:borne-M-infty} yields
\begin{eqnarray*}
\hspace*{-2mm} &&  |x|^{-p}   \PE_x [ M_\infty^p(\bepsilon{K}{},\lfloor T |\Phi_0|^{1+\beta}
  \rfloor) \1 \{ M_\infty(\bepsilon{K}{},T_{\Phi_0}) \geq (b/2) |\Phi_0| \}]
  \\
\hspace*{-2mm}&&\qquad  \leq (b/2)^{-(a-1)p} |x|^{-a p}
  \PE_x[M_\infty^{ap}(\bepsilon{K}{},T_{\Phi_0})] \leq C A(\bepsilon{K}{},ap)
  b^{-(a-1)p} |x|^{-a (1-\beta) p/2},
\end{eqnarray*}
where $A(\bepsilon{K}{},ap) \eqdef \sup_{x \in \Xset} \PE_x [|\bepsilon{K}{1}|^{ap}]$. Note that, by construction, $A(\bepsilon{K}{},ap) \leq K^{ap}$.
Similarly, Lemma \ref{lem:borne-M-infty} implies that
$ \PE_x [ M_\infty^p(\teps{K}{},T_{\Phi_0})] \leq C A(\teps{K},p)\times
T^{p/2} |x|^{\{p(1+\beta)/2\} \vee (1+\beta)}$, where
$A(\teps{K},p) \eqdef  \sup_{x \in \Xset} \PE_x [|\tepsilon{K}{1}|^{p}]$. Therefore, since $p \geq 1+\beta$,
$\sup_{|x|\geq 1} |x|^{-p} \PE_x [ M_\infty^p(\teps{K},T_{\Phi_0})] \leq C T^{p/2} A(\teps{K},p)$.
Combining the two last inequalities, we have
\begin{eqnarray*}
&& \sup_{|x|\geq 1} |x|^{-p} \PE_x[ M_\infty^p(\epsilon,T_{\Phi_0}) \1 \{
  M_\infty(\epsilon,T_{\Phi_0}) \geq b |\Phi_0| \} ]
  \nonumber
  \\[-8pt]
  \\[-8pt]
  \nonumber
  &&\qquad
   \leq C \{ K^{ap}
  b^{-p(a-1)} + A(\teps{K},p) \},
\end{eqnarray*}
which goes to 0 by setting $K \eqdef K(b) = \log(b)$.
\end{pf}

\begin{prop}
\label{prop:uniform-integrability}
Assume \textup{B1} and \textup{B2}.
Then, for all $T > 0$,
\begin{equation}
\label{eq:borne-Lp-Phi}
\sup_{x \in \Xset}   (1+|x|)^{-p}  \PE_x  \biggl[ \sup_{0 \leq k \leq \lfloor T |\Phi_0|^{1+\beta} \rfloor} |\Phi_k|^p  \biggr]  < \infty,
\end{equation}
\begin{eqnarray}\label{eq:uniforme-integrability-markov-chain}
&& \lim_{K \to \infty} \sup_{ |x| \geq 1} |x|^{-p} \PE_x  \biggl[ \sup_{0 \leq k
    \leq \lfloor T|\Phi_0|^{1+\beta}\rfloor} |\Phi_k|^p
    \nonumber
    \\[-8pt]
    \\[-8pt]
    \nonumber
&&  \hspace*{31mm} {}\times\1  \biggl\{ \sup_{0 \leq
      k \leq \lfloor T |\Phi_0|^{1+\beta}\rfloor} | \Phi_k | \geq K |\Phi_0|
   \biggr\}  \biggr] =  0.
\end{eqnarray}
\end{prop}

\begin{pf}
Set $T_{\Phi_0} = \lfloor T |\Phi_0|^{1+\beta} \rfloor$.  For all $r \geq 1$, applying Lemma~\ref{lem:ControleFluctuationPhi}(i)
with $K=J= \lfloor T |\Phi_0|^{1+\beta} \rfloor$ and $\kappa = |\Phi_0|$ yields
\begin{equation}
 \label{eq:majoPhi}
 \sup_{0 \leq k \leq T_{\Phi_0}} |\Phi_k|^r \leq C    \{ 1 +  |\Phi_0|^{r} + M_\infty^r(\epsilon,T_{\Phi_0})  \}
 \end{equation}
 for some constant $C$ depending upon $r,\beta,N(\beta,\Delta)$ and $T$.  The first assertion is
 then a consequence of Lemma~\ref{lem:borne-M-infty}.  Inequality \eqref{eq:majoPhi} applied with $r=1$ implies that there exist
 constants $a,b>0$ such that for all $|x| \geq 1$ and all large enough $K$,
 \[
  \biggl\{\sup_{0 \leq k \leq T_{\Phi_0}} |\Phi_k| \geq K |\Phi_0|  \biggr\}
  \subset  \{  M_\infty(\epsilon, T_{\Phi_0})  \geq (a K-b) |\Phi_0|
  \}\qquad \PP_x\mbox{-a.s.}
 \]
Hence, for large enough $K$ and an appropriately chosen constant $C$,
\begin{eqnarray*}
&& \sup_{|x| \geq 1} |x|^{-p} \PE_x  \biggl[ \sup_{0 \leq k \leq T_{\Phi_0}}  |\Phi_k|^p  \1  \biggl\{ \sup_{0 \leq k \leq T_{\Phi_0}} |\Phi_k| \geq K |\Phi_0|  \biggr\}  \biggr]
\\
&&\qquad  \leq  C \sup_{|x|\geq 1} \PP_x  [ M_\infty(\epsilon,T_{\Phi_0}) \geq (a K-b) |\Phi_0|  ]
\\
&&\qquad\quad{}   + C \sup_{|x|\geq 1} |x|^{-p}   \PE_x  [ M^p_\infty(\epsilon,T_{\Phi_0}) \1 \{ M_\infty(\epsilon,T_{\Phi_0}) \geq (a K-b)|\Phi_0| \} ].
\end{eqnarray*}
The proof of (\ref{eq:uniforme-integrability-markov-chain}) follows from Lemma \ref{lem:uniforme-integrability-martingale}.
\end{pf}

\begin{prop}
\label{theo:SuffCond-Drift-partII}
Assume \textup{B1} and
\textup{B2} and that there exist $T < \infty$ and $\rho
\in (0,1)$ such that
\begin{equation}
\label{eq:key-condition-stopping-time}
\limsup_{|x| \to \infty} \PP_x  ( \sigma > \tau  ) = 0,\qquad\mbox{with }   \sigma \eqdef \inf \{ k \geq 0, |\Phi_k| < \rho |\Phi_0|  \},
\end{equation}
where  $\tau \eqdef \sigma \wedge  \lceil T|\Phi_0|^{1+\beta} \rceil$.  It then follows that \textup{(a)} there
exists $M$ such that $\sup_{|x| \geq M} |x|^{-p} \PE_x  [ |\Phi_{\tau}|^p  ] < 1$   and

\textup{(b)} $\PE_{x}  [ \sum_{k=0}^{\tau -1} |\Phi_k|^p  ] \leq C   |x|^{p+1+\beta}$.
\end{prop}

\begin{pf}
Set $T_{\Phi_0} = \lceil T |\Phi_0|^{1+\beta} \rceil$. For any $K \geq
0$,
\begin{eqnarray}
\label{eq:decomposition}
\nonumber && |x|^{-p} \PE_x  [  | \Phi_{\tau}  |^p ]
\\
&&\qquad =  |x|^{-p}  \PE_x
 [ \1 {\{ \tau = \sigma\}}  | \Phi_{\tau} |^p ] +
|x|^{-p} \PE_x  [ \1 {\{ \sigma > T_{\Phi_0} \} } | \Phi_{T_{\Phi_0}} |^p ]
\\
\nonumber
&& \qquad \leq {\rho}^{p} + |x|^{-p} \PE_x  [  | \Phi_{T_{\Phi_0}}  |^p \1
  \{ | \Phi_{T_{\Phi_0}}| \geq K |\Phi_0| \}  ] + K^p \PP_{x} [ \sigma >
T_{\Phi_0}].
\end{eqnarray}
By Proposition~\ref{prop:uniform-integrability}, one may choose $K$
sufficiently large so that
\begin{equation}
\label{eq:choix-de-K}
\sup_{|x|\geq 1} |x|^{-p} \PE_x  [   | \Phi_{ T_{\Phi_0}}  |^p \1 \{ |\Phi_{T_{\Phi_0}}| \geq K |\Phi_0| \}  ] < 1 - {\rho}^p.
\end{equation}
Since $\limsup_{|x| \to \infty} \PP_x [\sigma > T_{\Phi_0}]= 0$, the proof of (a) follows.
Since $\tau \leq T_{\Phi_0}$,  (b) follows from  \eqref{eq:borne-Lp-Phi} and the bound
$\PE_{x}  [ \sum_{k=0}^{\tau -1} |\Phi_k|^p  ] \leq C   T |x|^{1+\beta}\times\break   \PE_x [ \sup_{1 \leq k \leq T_{\Phi_0}} |\Phi_k|^p ]$.
\end{pf}

The following elementary proposition relates the stability of the fluid limit
model to the condition \eqref{eq:key-condition-stopping-time} on the stopping
time $\sigma$. We introduce the polygonal process that agrees with $\Phi_k/r$
at the points $t = k r^{-(1+\alpha)}$ and is defined by linear interpolation
 \begin{eqnarray}
   \label{eq:interpolated-continuous-process}
   \tilde{\eta}_r^{\alpha}({t};{x})&=&  r^{-1} \sum_{k \geq 0}   \{
   ( k+1-t r^{1+\alpha} ) \Phi_{k} +  (t r^{1+\alpha} -k  ) \Phi_{k+1}  \}
   \nonumber
   \\[-16pt]
   \\
   \nonumber
   &&\hspace*{10.5mm}
{}\times
   \1 \{  k  \leq t r^{1+\alpha} < (k+1)\}.
 \end{eqnarray}
 Denote by $\tilde{\mathbb{Q}}^{\alpha}_{r;x}$ the image probability on
 $\continuousf{\rset^+,\Xset}$ of $\PP_{rx}$ by $\tilde\eta_r^{\alpha}({t};{x})$.  The introduction of this process allows for an easier
 characterization of the open and closed sets of
 $\continuousf{[0,T],\Xset}$ equipped with the uniform topology, than
 the open and closed sets of $\mathsf{D}([0,T],\Xset)$ equipped with the
 Skorokhod topology.  For any sequences $\{r_n\}_n \subset \rset^+$ such that
 $r_n \to +\infty$ and $\{x_n \} \subset \Xset$ such that $x_n \to x$, the
 family of probability measures $\{\tilde{\mathbb{Q}}_{r_n;x_n}^{\alpha}\}$ is tight and
 converges weakly to $\mathbb{Q}^{\alpha}_{x}$, the weak limit of the sequence
 $\{\mathbb{Q}_{r_n;x_n}^{\alpha}\}_{n \in \nset}$.  This can be proved following the
 same lines as in the proof of Theorem~\ref{theo:weak-compacity} [see,
 e.g., \citet{billingsley1999}, Theorem~7.3].  Details are
 omitted.

\begin{prop}
\label{prop:stability-fluid-limit}
Assume \textup{B1, B2} and that the $\beta$-fluid limit model $\{\mathbb{Q}^\beta_x, x \in \Xset \}$ is stable.
Then,  \eqref{eq:key-condition-stopping-time} is satisfied.
\end{prop}

\begin{pf}
Let $\{y_n\} \subset \Xset$ be any sequence of initial states with $|y_n| \to \infty$ as $n \to \infty$. Set $r_n \eqdef |y_n|$ and $x_n \eqdef y_n/|y_n|$. One may extract
a subsequence $\{x_{n_j} \} \subseteq \{ x_n \}$ such that $\lim_{j \to \infty} x_{n_j} =x$ for some $x$, $|x|=1$.
By Theorem \ref{theo:weak-compacity}, there exist  subsequences $\{r_{m_j} \}\subseteq \{r_{n_j}\}$ and $\{ x_{m_j} \} \subseteq \{ x_n \}$ and a $\beta$-fluid limit
$\mathbb{Q}^\beta_x$ such that $\tilde{\mathbb{Q}}_{r_{m_j}; x_{m_j}}^\beta \Rightarrow \mathbb{Q}^\beta_x$. By construction,
\begin{eqnarray*}
\PP_{r_{m_j} x_{m_j}}  (\sigma > \tau  )
&\leq& \PP_{r_{m_j} x_{m_j}}  \biggl( \inf_{0 \leq t \leq T} | \tilde\eta_{r_{m_j}}^{\beta}({t};{x_{m_j}})| \geq \rho  \biggr)
\\
&=&  \tilde{\mathbb{Q}}_{r_{m_j}; x_{m_j}}^\beta  \biggl( \eta \in \mathsf{C}(\rset^+,\Xset)\dvtx  \inf_{0 \leq t \leq T} |\eta(t)| \geq \rho  \biggr).
\end{eqnarray*}
By the Portmanteau theorem, since the set $\{\eta \in \mathsf{C}(\rset^+,\Xset),   \inf_{[0,T]} |\eta| \geq \rho\}$ is
closed, we have
\[
\limsup_{j \to \infty} \tilde{\mathbb{Q}}_{r_{m_j}; x_{m_j}}^\beta \biggl( \inf_{0 \leq t \leq T} |\eta(t)| \geq \rho
 \biggr) \leq \mathbb{Q}^\beta_x  \biggl( \inf_{0 \leq t \leq T} |\eta(t)| \geq \rho  \biggr) = 0.
\]
Because $\{y_n\}$ is an arbitrary sequence, this relation implies \eqref{eq:key-condition-stopping-time}.
\end{pf}

\begin{pf*}{Proof of Theorem \ref{theo:(f,r)-ergodicity}}
  This follows immediately from Theorem \ref{theo:state-dependent-drift}, using
  Propositions \ref{theo:SuffCond-Drift-partII} and
  \ref{prop:stability-fluid-limit}.
\end{pf*}

\subsection[Proof of Proposition~1.5]{Proof of Proposition~\textup{\protect\ref{prop:characterization-fluid-limit}}}
\label{sec:proof:prop:characterization-fluid-limit}
In this proof, we see the $\beta$-fluid limit $\mathbb{Q}^\beta_x$ as the weak limit
of $\tilde{\mathbb{Q}}_{r_n;x_n}^{\beta}$ for some sequences $\{r_n\} \subset \rset_+$
and $\{x_n\} \subset \Xset$ satisfying $\lim_{n \to \infty} r_n = \infty$ and
$\lim_{n \to \infty} x_n = x$. Fix $s,t$ such that $s<t$.  We prove that
 \begin{eqnarray}
   \label{eq:ProbaCondNulle}
&&   \mathbb{Q}^\beta_x  \biggl( \mathsf{A}(s,t) \cap  \biggl\{ \eta \in \mathsf{C}([s,t],\Xset) \dvtx
\nonumber
\\[-8pt]
\\[-8pt]
\nonumber
&&\hspace*{24.2mm}   \sup_{s \leq u
    \leq t}  \bigg| \eta(u) -\eta(s) - \int_s^u h \circ \eta(y)\,dy  \bigg| >0  \biggr\}  \biggr) =0.
 \end{eqnarray}
Let $\mathsf{U}$ be an  open set such that $\bar{\mathsf{U}} \subseteq \openset$, where $\bar{\mathsf{U}}$ denotes the
closure of the set $\mathsf{U}$.  For any $\delta > 0$, $M>0$ and $m >0$, $ s \leq u < w \leq t$, define
\begin{eqnarray}
\label{eq:definition-Adelta}
\mathsf{A}^{\mathsf{U}}_{\delta,m,M}(u,w) &\eqdef&  \biggl\{ \eta \in \mathsf{C}([s,t],\Xset),  \eta([u,w]) \subset \mathsf{U} \cap \mathsf{C}_{m, M},
\nonumber
\\[-8pt]
\\[-8pt]
\nonumber
&&\phantom{\bigg\{}
 \sup_{u \leq v \leq w}  \bigg| \eta(v) - \eta(u) - \int_u^v h \circ \eta(x)\,dx  \bigg| > \delta  \biggr\},
\end{eqnarray}
where $\mathsf{C}_{m,M} \eqdef  \{ x \in \Xset, m \leq |x| \leq M  \}$.
Since $\delta$, $m$, $M$, $\mathsf{U}$, $u$ and $w$ are arbitrary,
(\ref{eq:ProbaCondNulle}) holds whenever $\mathbb{Q}^\beta_x [\mathsf{A}^{\mathsf{U}}_{\delta,m,M}(u,w)  ] = 0$.
By the Portmanteau theorem, since the set $\mathsf{A}_{\delta,m,M}^{\mathsf{U}}(u,w)$ is open in the uniform topology,
\[
\mathbb{Q}^\beta_x [ \mathsf{A}_{\delta,m,M}^{\mathsf{U}}(u,w)  ] \leq \liminf_{n \to \infty} \tilde{\mathbb{Q}}_{r_n;x_n}^{\beta}  [ \mathsf{A}_{\delta,m,M}^{\mathsf{U}}(u,w)  ]
\]
and the property will follow if we can prove that the right-hand side of\vadjust{\goodbreak} the previous inequality is null.
To that end, we write
\begin{eqnarray*}
  && \tilde{\eta}_{r_n}^{\beta}({v};{x_n}) - \tilde\eta_{r_n}^{\beta}({u};{x_n}) - \int_u^v h \circ \tilde\eta_{r_n}^{\beta}({y};{x_n}) \,dy
  \\
  &&\qquad = \tilde\eta_{r_n}^{\beta}({v};{x_n}) -
  \tilde\eta_{r_n}^{\beta}\bigl({\lfloor v r_n^{1+\beta}\rfloor
  r_n^{-(1+\beta)}};{x_n}\bigr)
\\
  &&\quad\qquad{}+ \tilde\eta_{r_n}^{\beta}\bigl({\lfloor u r_n^{1+\beta}\rfloor r_n^{-(1+\beta)}};{x_n}\bigr)
  - \tilde\eta_{r_n}^{\beta}({u};{x_n})
  \\
  &&\qquad\quad{}  +r_n^{-1} \sum_{k= \lfloor u r_n^{1+\beta } \rfloor}^{\lfloor v r_n^{1+\beta } \rfloor -1} \{ \Phi_{k+1} - \Phi_k \} - \int_u^v h \circ \tilde\eta_{r_n}^{\beta}({t};{x_n}) \,dt
  \\
  &&\qquad \leq 2 \chi_1 + \chi_2 + \chi_3 + 2 r_n^{-1} M_\infty(\epsilon,\lfloor t r_n^{1+\beta } \rfloor),
\end{eqnarray*}
where we have defined
\begin{eqnarray*}
 \chi_1 & \eqdef&  \sup_{u \leq v \leq w}  \biggl\{  \big| \tilde\eta_{r_n}^\beta (v;{x_n}) -
 \tilde{\eta}_{r_n}^{\beta}\bigl(\lfloor v r_n^{1+\beta}\rfloor r_n^{-(1+\beta)};{x_n}\bigr)  \big|
\\
&&\hspace*{18mm}
{} +  \bigg| \int_{\lfloor v r_n^{1+\beta}\rfloor r_n^{-(1+\beta)}}^v h \circ \tilde\eta_{r_n}^{\beta}({t};{x_n}) \,dt  \bigg|  \biggr\},
 \\
 \chi_2 &=&  \sum_{j= \lfloor u r_n^{1+\beta} \rfloor}^{\lfloor w r_n^{1+\beta} \rfloor -1}
 \big|  r_n^{-1} \Delta\bigl( r_n \tilde\eta_{r_n}^{\beta}\bigl(j r_n^{-(1+\beta)}\bigr);{x_n}\bigr)
 - r_n^{-(1+\beta)} h \bigl( \tilde\eta_{r_n}^{\beta}\bigl(j r_n^{-(1+\beta)}\bigr){x_n} \bigr)  \big|,
 \\
 \chi_3 &=&  \sum_{j= \lfloor u r_n^{1+\beta} \rfloor}^{\lfloor w r_n^{1+\beta} \rfloor -1}
 \bigg| r_n^{-(1+\beta)} h\bigl( \tilde\eta_{r_n}^{\beta}\bigl(j r_n^{-(1+\beta)};{x_n}\bigr) \bigr)
 - \int_{j r_n^{-(1+\beta)} }^{(j+1) r_n^{-(1+\beta)}} h \circ \tilde\eta_{r_n}^{\beta}({t};{x_n})\,  dt  \bigg|.
\end{eqnarray*}
Denote by  $\omega_{m,M,\mathsf{U}}$ the modulus of continuity of $h$ on
$\mathsf{U} \cap \mathsf{C}_{m,M}$. Since $h$ is continuous on $\mathsf{U}$, $\lim_{\lambda \to 0}
\omega_{m,M,\mathsf{U}}(\lambda)=0$. On the event $ \{ \tilde\eta_{r_n}^{\beta}({t};{x_n}) \in \mathsf{U} \cap \mathsf{C}_{m,M}  \}$,
\begin{eqnarray*}
\chi_1 &\leq&  r_n^{-1}  \biggl( 1+ \sup_{|x| \geq m} |h(x)|  \biggr) \sup_{1 \leq j \leq \lfloor t r_n^{1+\beta} \rfloor } |\Phi_{j+1} - \Phi_j |,
\\
\chi_2 &\leq& (t-s+1)  m^{-\beta} \sup_{\{ x \in \mathsf{U}, |x| \geq m \}}  \big|   r_n^\beta |x|^\beta \Delta(r_n x) - \Delta_\infty(x)    \big|
\end{eqnarray*}
and, for any $\lambda > 0$,
\[
\chi_3  \leq (t-s+1)    \biggl( \omega_{m,M,\mathsf{U}}(\lambda) + \sup_{|x| \geq m} |h(x)|   \1
\biggl\{ \sup_{1 \leq j \leq \lfloor t r_n^{1+\beta} \rfloor } | \Phi_{j+1} - \Phi_j| \geq \lambda r_n  \biggr\} \biggr).
\]
By Lemma \ref{lem:ControleFluctuationPhi}, for any $\delta > 0$, $\lim_{n \to \infty} \PP_{r_n x_n}  ( \sup_{1 \leq j \leq \lfloor t r_n^{1+\beta} \rfloor } |\Phi_{j+1} - \Phi_j |  \geq \delta r_n  )= 0$.
On the other hand, $\lim_{n \to \infty} \sup_{\{ x \in \mathsf{U}, |x| \geq m r_n\}}  |   |x|^\beta \Delta(x) - \Delta_\infty(x)    |= 0$.
Therefore, for any $\delta > 0$, one may choose $\lambda$ small enough so that
\[
\lim_{n \to \infty} \PP_{r_n x_n}  \bigl( \tilde\eta_{r_n}^{\beta}({t};{x_n}) \in \mathsf{U} \cap \mathsf{C}_{m,M}, (2 \chi_1 + \chi_2 + \chi_3) \geq \delta   \bigr) = 0.
\]
The proof follows from Lemma \ref{lem:borne-M-infty}.

\subsection[Proof of Theorem~1.6]{Proof of Theorem~\textup{\protect\ref{theo:fluidlimit-Lipsch-NeighLimit}}}
\label{sec:ProofCoro1}
We preface the proof by a lemma showing that the fluid limits are uniformly
bounded.

\begin{lem}
\label{lem:ControleEcartLF}
Assume \textup{B1} and \textup{B2}.
\begin{longlist}[(ii)]
\item[(i)] For any $T>0$ and $\rho>0$, there exists $\delta>0$ such that, for any
  $\beta$-fluid limit~$\mathbb{Q}^\beta_x$,
\begin{equation}
\label{eq:control-modulus-of-continuity}
\mathbb{Q}^\beta_x \biggl( \eta \in \mathsf{C}(\rset^+,\Xset), \sup_{0 \leq t \leq u \leq t +\delta \leq T} | \eta(u) - \eta(t)| \leq \rho  \biggr) =1.
\end{equation}

\item[(ii)] For any $T>0$, there exists $K>0$ such that, for any $\beta$-fluid limit~$\mathbb{Q}^\beta_x$,
\begin{equation}
\label{eq:control-norm}
\mathbb{Q}^\beta_x \biggl( \eta \in \mathsf{C}(\rset^+,\Xset),   \sup_{0 \leq t \leq T} |\eta(t) - \eta(0)| \geq K  \biggr) =0.
\end{equation}
\end{longlist}
\end{lem}

\begin{pf}
(i)  Let $\{r_n \} \subset \rset_+$ and $\{ x_n \} \subset \Xset$
be two sequences such that\break $\lim_{n \to \infty} r_n = +\infty$, $\lim_{n \to \infty} x_n = x$ and $\mathbb{Q}_{r_n; x_n}^{\beta} \Rightarrow \mathbb{Q}^\beta_x$.
 By the Portmanteau theorem, since the set $\{\eta \in \mathsf{C}(\rset^+,\Xset),   \sup_{0 \leq t \leq u \leq t+\delta \leq T}
  |\eta(u) -\eta(t)| \leq \rho \}$ is closed, it follows that
  \begin{eqnarray*}
    && \mathbb{Q}^\beta_x \biggl( \eta \in \mathsf{C}(\rset^+,\Xset),  \sup_{0 \leq t \leq u \leq t+\delta \leq T} | \eta(u) - \eta(t)| \leq \rho  \biggr)
    \\
    &&\qquad \geq \limsup_n \tilde{\mathbb{Q}}_{r_n;x_n}^{\beta} \biggl( \eta \in
      \mathsf{C}(\rset^+,\Xset),  \sup_{0 \leq t \leq u \leq t+\delta \leq T} | \eta(u) - \eta(t)| \leq \rho  \biggr).
  \end{eqnarray*}
By definition of the process $\tilde\eta_{r_n}^{\beta}({\cdot};{x_n})$,
\begin{eqnarray*}
&&  \tilde{\mathbb{Q}}_{r_n;x_n}^{\beta} \biggl(\eta \in \mathsf{C}(\rset^+,\Xset), \sup_{0 \leq t \leq u \leq t+\delta \leq T} | \eta(u) - \eta(t)| > \rho  \biggr)
\\
&&\qquad  \leq \PP_{r_n x_n}  \biggl(\sup_{0 \leq k < k+j \leq T r_n^{1+\beta}, 0 \leq j
      \leq \delta r_n^{1+\beta}} | \Phi_{k+j} - \Phi_k| > \rho r_n \biggr)
\end{eqnarray*}
and the proof follows from
Lemma~\ref{lem:ControleFluctuationPhi}(iii).

(ii) The proof follows from (i) by considering  the decomposition
\[
\sup_{0 \leq t \leq T} |\eta(t) - \eta(0)|  \leq \sum_{q=0}^{\lfloor T/\delta \rfloor} \sup_{q \delta \leq u \leq (q+1) \delta} |\eta(u) -\eta(q \delta)|.
\]\upqed
\end{pf}

\begin{pf*}{Proof of Theorem~\ref{theo:fluidlimit-Lipsch-NeighLimit}}
Under the stated assumptions, $\mu([0,T_x];)$ is a compact subset of $\openset$.
Since $\openset$ is open, there exists $\rho>0$ such that
\[
 \{ y \in \Xset, d  (y, \mu([0,T_x];x)  ) \leq 2 \rho  \} \subset \openset,
\]
where, for $x \in \Xset$ and $A \subset \Xset$, $d(x,A)$ is the distance from
$x$ to the set $A$.  By Lemma~\ref{lem:ControleEcartLF}(i), there
exists $\delta > 0$ such that
\[
\mathbb{Q}^\beta_x \biggl(\eta \in \mathsf{C}(\rset_+,\Xset), \sup_{0 \leq t \leq u \leq t+\delta \leq T_x} |\eta(u) - \eta(t)| \leq \rho  \biggr) =1.
\]
Since $\mathbb{Q}^\beta_x (\eta \in \mathsf{C}(\rset_+,\Xset),  \eta(0) = x = \mu(0;x)
)=1$, we have
\[
\mathbb{Q}^\beta_x \bigl(\eta \in \mathsf{C}(\rset_+,\Xset),  \eta([0, \delta]) \subset \openset  \bigr)=1.
\]
By Proposition~\ref{prop:characterization-fluid-limit}, this yields $\mathbb{Q}^\beta_x = \delta_{\mu(\cdot;x)}$ on
$\mathsf{C}([0,\delta],\Xset)$.  By repeated application of
Lemma~\ref{lem:ControleEcartLF}(i), it is readily proved by induction
that $\mathbb{Q}^\beta_x = \delta_{\mu(\cdot;x)}$ on $\mathsf{C} ([(q-1) \delta,q\delta] \cap[0,T_x],\Xset)$ for any integer $q \geq 1$.
\end{pf*}

\subsection[Proof of Theorem~1.9]{Proof of Theorem~\textup{\protect\ref{theo:stability-fluid-limit-deterministic-nonsmooth}}}
\label{sec:proof:theo:stability-fluid-limit-deterministic-nonsmooth}
Let $x$ be such that $|x|=1$.  By Lemma~\ref{lem:ControleEcartLF}, there exists
$K$ depending on $T_0$ such that $\mathbb{Q}^\beta_x  (\eta\dvtx
\sup_{[0,T_0]} |\eta(\cdot)| \leq K  )=1$ for any $\beta$-fluid
limit $\mathbb{Q}^\beta_x$.  Set $T = T_0 +T_K$, where $T_0$ and $T_K$ are defined by
\eqref{eq:IntersectCGamma} and \eqref{eq:DefiTK}, respectively.

By definition, for any set $\compactset$, $\compactset \subset
\Omega_\compactset$; therefore, there exists an increasing sequence $\{
\compactset_n \}$ of compact subsets of $\openset$ such that $\compactset_n
\subsetneq \compactset_{n+1}$ and $\openset= \bigcup_n \Omega_{\compactset_n}$
(note that $\Omega_{\compactset_n} \subseteq \Omega_{\compactset_{n+1}}$). This
implies that
\begin{eqnarray*}
&&   \mathbb{Q}^\beta_x \biggl(\eta\dvtx  \inf_{[0,T]} |\eta(\cdot)| >
    \rho_K  \biggr)
\\
&&\qquad =  \mathbb{Q}^\beta_x \biggl(\eta\dvtx
    \inf_{[0,T]} |\eta(\cdot)| > \rho_K, \eta([0,T_0]) \cap \openset \neq \varnothing  \biggr)
    \\
&&\qquad =  \lim \uparrow_n \mathbb{Q}^\beta_x \biggl(\eta\dvtx  \inf_{[0,T]}
    |\eta(\cdot)| > \rho_K, \eta([0,T_0]) \cap
    \Omega_{\compactset_n} \neq \varnothing  \biggr).
\end{eqnarray*}
$\lim\uparrow_n$ stands for a limit that converges monotonically from below.
We prove that for any $n$, the term in the right-hand side is zero. To that
end we start by proving that for any compact set $ \compactset \subset
\openset$ and any real numbers $ 0 \leq q \leq T_0$,
\begin{eqnarray}  \label{eq:NestedReal}
\nonumber   && \mathbb{Q}^\beta_x \biggl(\eta\dvtx  \inf_{[0,T]} |\eta(\cdot)| > \rho_K, \eta(q) \in \Omega_{\compactset}  \biggr)
\\
  &&\qquad = \mathbb{Q}^\beta_x \biggl(\eta\dvtx  \inf_{[0,T]} |\eta(\cdot)| > \rho_K, \eta( q+\cdot) =
    \mu({\cdot};\eta(q))
\\
\nonumber &&\hspace*{129pt}   \mbox{ on }   \bigl[0,T_{\eta(q)} \bigr], \eta(q) \in
    \Omega_{\compactset}  \biggr).
\end{eqnarray}
We will then establish that
\begin{eqnarray}
  \label{eq:IntervalQ}
&& \mathbb{Q}^\beta_x \biggl(\eta\dvtx  \inf_{[0,T]} |\eta(\cdot)| > \rho_K, \eta( q+\cdot) =
    \mu({\cdot};\eta(q))
    \nonumber
    \\[-8pt]
    \\[-8pt]
    \nonumber
    &&\hspace*{95.5pt}    \mbox{ on }   \bigl[0,T_{\eta(q)} \bigr], \eta(q) \in
    \Omega_{\compactset}  \biggr) =0.
\end{eqnarray}
Since
$\mathbb{Q}^\beta_x(\mathsf{C}(\rset^+,\Xset))=1$, \eqref{eq:NestedReal} and \eqref{eq:IntervalQ} imply that
\begin{eqnarray*}
&&  \mathbb{Q}^\beta_x \biggl( \eta\dvtx  \inf_{[0,T]} |\eta(\cdot)| > \rho_K, \eta([0,T_0]) \cap \Omega_{\compactset_n} \neq
    \varnothing  \biggr)
    \\
&&\qquad   \leq \sum_{q \in \mathcal{Q}} \mathbb{Q}^\beta_x \biggl(\eta\dvtx  \inf_{[0,T]} |\eta(\cdot)| > \rho_K,
    \eta(q+\cdot) = \mu({\cdot};\eta(q))
\\
&&\hspace*{141.6pt}      \mbox{ on } \bigl[0, T_{\eta(q)} \bigr],
    \eta(q) \in \Omega_{\compactset'_n}  \biggr) = 0,
\end{eqnarray*}
where $\compactset'_n \supset \compactset_n$ is a compact set of $\openset$ and $\mathcal{Q} \subset [0,T_0]$ is a denumerable dense  set.
This  concludes the proof.

We now turn to the proof of (\ref{eq:NestedReal}) and (\ref{eq:IntervalQ}).
Since $\Omega_{\compactset}$ is a compact set of $\openset$, there exists
$\varepsilon>0$ (depending on $\compactset$) such that $\{y \in \Xset,
d ( y,  \Omega_{\compactset}  ) \leq 2 \varepsilon \} \subsetneq
\openset$.  By Lemma~\ref{lem:ControleEcartLF}, one may choose $\delta > 0$
small enough (depending on $T$ and $\varepsilon$) so that
\[
\mathbb{Q}^\beta_x \biggl( \eta \in \mathsf{C}(\rset^+,\Xset)\dvtx  \sup_{0 \leq t \leq u \leq t
    +\delta \leq T} | \eta(u) - \eta(t)| \leq \varepsilon  \biggr) =1.
\]
Therefore, for any compact set $\compactset \subset \openset$ and $q \in \mathcal{Q}$,
\begin{eqnarray*}
&& \mathbb{Q}^\beta_x \biggl(\eta \dvtx  \inf_{[0,T]} |\eta(\cdot)| > \rho_K,  \eta(q) \in \Omega_{\compactset}  \biggr)
 \\
&&\qquad   = \mathbb{Q}^\beta_x \biggl(\eta\dvtx  \inf_{[0,T]} |\eta(\cdot)| > \rho_K,
  \eta(q) \in \Omega_{\compactset},   \sup_{0 \leq t \leq u \leq t + \delta \leq T} |\eta(u) -\eta(t)| \leq \varepsilon  \biggr)
  \\
&&\qquad   = \mathbb{Q}^\beta_x \biggl(\eta\dvtx
\inf_{[0,T]} |\eta(\cdot)| > \rho_K,  \eta(q) \in \Omega_{\compactset},   \eta\bigl([q, (q+\delta)\wedge T ]\bigr) \subset \openset \biggr).
\end{eqnarray*}
By Proposition~\ref{prop:characterization-fluid-limit}, on the set $\mathsf{A}(q,q+\delta)$,
$\eta(q+ \cdot)= \mu({\cdot};\eta(q))$ on $[0,\delta \wedge T_{\eta(q)}]$, $\mathbb{Q}^\beta_x$-a.s.
Hence,
\begin{eqnarray*}
&&  \mathbb{Q}^\beta_x \biggl(\eta \dvtx  \inf_{[0,T]} |\eta(\cdot)| > \rho_K, \eta(q) \in \Omega_{\compactset}  \biggr)
\\
&&\qquad  = \mathbb{Q}^\beta_x \biggl(\eta \dvtx  \inf_{[0,T]} |\eta(\cdot)| > \rho_K, \eta(q) \in \Omega_{\compactset}, \eta(q+\cdot)
    = \mu({\cdot};\eta(q)),
\\
&&\hspace*{79.2mm}\hspace*{-8.2pt}      \mbox{ on }  \bigl[0,\delta \wedge T_{\eta(q)}\bigr]  \biggr).
\end{eqnarray*}
By repeated application of Proposition~\ref{prop:characterization-fluid-limit}, for any integer $l>0$,
\begin{eqnarray*}
&&  \mathbb{Q}^\beta_x \biggl(\eta \dvtx  \inf_{[0,T]} |\eta(\cdot)| > \rho_K, \eta(q) \in \Omega_{\compactset}  \biggr)
\\
&&\qquad  = \mathbb{Q}^\beta_x \biggl(\eta \dvtx  \inf_{[0,T]} |\eta(\cdot)| > \rho_K, \eta(q) \in
    \Omega_{\compactset}, \eta(q+\cdot) = \mu({\cdot};\eta(q)),
\\
&&\hspace*{78.1mm} \hspace*{-9.3pt}      \mbox{ on }  \bigl[0,l \delta \wedge
    T_{\eta(q)}\bigr]  \biggr),
\end{eqnarray*}
which concludes the proof of \eqref{eq:NestedReal}.
\begin{eqnarray*}
&&  \mathbb{Q}^\beta_x \biggl(\eta\dvtx  \inf_{[0,T]} |\eta(\cdot)| > \rho_K, \eta(q+ \cdot) =
\mu({\cdot};\eta(q))
     \mbox{ on } \bigl[0, T_{\eta(q)}\bigr], \eta(q) \in \Omega_{\compactset} \biggr)
     \\
&&\qquad \leq \mathbb{Q}^\beta_x \biggl(\eta\dvtx  \inf_{[0,T]} |\eta(\cdot)| > \rho_K, \inf_{[0, T_0+T_K]} |\eta| \leq \rho_K  \biggr) =0
\end{eqnarray*}
since $T = T_0 + T_K$, which concludes the proof of (\ref{eq:IntervalQ}).

\section[Proofs for  Section~2]{Proofs for  Section~\textup{\protect\ref{sec:The ODE method for the Metropolis-Hastings algorithm}}}

\subsection[Proofs of  Section~2.1]{Proofs of  Section~\textup{\protect\ref{subsec:Super-exponential target densities}}}
\label{proof:sec:superexp}

\mbox{}
\begin{pf*}{Proof of Proposition \ref{prop:super-exponential}}
Define
\begin{equation}
\label{eq:limiting-field-super-exponential}
\bar{\Delta}(x) \eqdef -\! \int_{\mathsf{R}_x} y  q(y) \lleb(dy).
\end{equation}
Introduce, for any $\delta > 0$, the \mbox{$\delta$-}zone $\mathsf{C}_{x}(\delta)$
around $\mathsf{C}_{x}$,
\begin{equation}
\label{eq:delta-radial-zone}
\mathsf{C}_{x}(\delta) \eqdef  \{ y + s n(y), y \in \mathsf{C}_{x}, -\delta \leq s \leq \delta  \}.
\end{equation}
By \citeauthor{jarnerhansen2000} [(\citeyear{jarnerhansen2000}), Theorem~4.1], we may bound the measure of the
\mbox{$\delta$-}zone's intersection with the ball $\mathsf{B}(0,K)$, for any $K > 0$
and all $|x|$ large enough,
\[
\lleb \bigl( \mathsf{C}_{x}(\delta) \cap \mathsf{B}(0,K)  \bigr)
\leq \delta  \biggl( \frac{|x|+K}{|x|-K} \biggr)^{d-1} \frac{\lleb \{ \mathsf{B}(0,3K) \}}{K},
\]
where the $x$-dependent term tends to $1$ as $|x|$ tends to infinity. From
this, it follows, using the fact that $\int |y| q(y) \lleb(dy) < \infty$, that for any $K > 0$
and $\epsilon > 0$, there exists $\delta > 0$ such that
\begin{equation}
\label{eq:bound}
\limsup_{|x| \to \infty} \int_{\mathsf{E}_x(\delta,K)}  |y| q(y) \lleb(dy) < \epsilon,
\end{equation}
where $\mathsf{E}_x(\delta,K) \eqdef  \mathsf{C}_{x}(\delta) \cap
  \mathsf{B}(0,K)$.  For arbitrary, but fixed, $\epsilon > 0$, choose
$K> 0$ such that $\int_{\mathsf{B}^c(0,K)} |y| q(y) \lleb(dy) \leq \epsilon$.
Then choose $\delta > 0$ such that \eqref{eq:bound} holds.  By construction,
for $y \in \mathsf{R}_x$, $\pi(x+y)/\pi(x) \leq 1$ and \eqref{eq:bound}
implies that
\begin{eqnarray}
\label{eq:premier-bout}
\limsup_{|x| \to \infty} \int_{\mathsf{R}_x \cap \mathsf{E}_x(\delta,K)} |y| \frac{\pi(x+y)}{\pi(x)} q(y) \lleb(dy) &\leq& \epsilon,
\\
\label{eq:second-bout}
\limsup_{|x| \to \infty} \int_{\mathsf{R}_x \cap \mathsf{B}^c(0,K)} |y| \frac{\pi(x+y)}{\pi(x)} q(y) \lleb(dy) &\leq& \epsilon.
\end{eqnarray}
From \eqref{eq:radial-small}, for $y \in \mathsf{R}_x$ such that $y$ has
radial distance at least $\delta$ to $\mathsf{C}_{x}$, the acceptance
probability satisfies $\pi(x+y)/\pi(x) \leq \epsilon /K$ for all $|x|$
sufficiently large [see~\citet{jarnerhansen2000}, page~351] and \eqref{eq:bound}
shows that
\begin{equation}
\label{eq:troisieme-bout}
\limsup_{|x| \to \infty} \int_{\mathsf{R}_x \cap \mathsf{E}_x^c(\delta,K)   \cap \mathsf{B}(0,K)} |y| \frac{\pi(x+y)}{\pi(x)} q(y) \lleb(dy) \leq \epsilon.
\end{equation}
By combining \eqref{eq:DefinitionDelta},
\eqref{eq:limiting-field-super-exponential}, \eqref{eq:premier-bout},
\eqref{eq:second-bout} and \eqref{eq:troisieme-bout}, $\limsup_{|x| \to \infty}
|\Delta(x) - \bar{\Delta}(x)| \leq 3 \epsilon$ and since $\epsilon$ is
arbitrary, $\lim_{|x|\to \infty} |\Delta(x) - \bar{\Delta}(x)| =0$.
\end{pf*}

\begin{pf*}{Proof of Proposition \protect\ref{prop:DeltaInftyExpReg}}
Set $z= (z_1, \dots, z_d) \eqdef \Sigma^{-1/2} y$ and $v = n(\Sigma^{1/2} u)$.
Then,
\begin{eqnarray*}
\int_{\{y, y'u \geq 0 \}}   y q (y) \lleb(dy)  &=&
\Sigma^{1/2} \int_{\{z, v'z \geq 0 \}}   z q_0(z) \lleb(dz)
\\
&=&  \Sigma^{1/2} v \int_{\Xset} z_1 \1_{\{z_1 \geq 0\}} q_0(z)\,dz.
\end{eqnarray*}
The proof follows.
\end{pf*}

\subsection[Proof of Lemma 2.9]{Proof of Lemma \textup{\protect\ref{lem:champ-de-trefle}}}
\label{sec:proof:MCMC-Ex2}
Let $\delta$ and $M$ be constants to be specified later. Write $\Delta(x) -
\Delta_\infty(x) \eqdef \sum_{i=1}^4 A_i(\delta,M,x)$, where
\begin{eqnarray*}
   A_1(\delta,M,x) &\eqdef& \int_{\{y, |y| \leq M,  | y' \Gamma_2^{-1} x  | \geq \delta |x|\}}  \frac{\pi(x +y) }{ \pi(x)} \1_{R_{\infty,  x}}(y) \  y q(y) \lleb(dy),
   \\
   A_2(\delta,M,x) &\eqdef& \int_{\{y, |y| \leq M,   | y' \Gamma_2^{-1} x  | \geq \delta |x| \} }
   \biggl( \frac{\pi(x +y) }{ \pi(x)}  -1  \biggr)\bigl( \1_{R_{x}}(y) - \1_{R_{\infty,x}}(y)  \bigr)
   \\
   &&\hspace*{33.3mm} {}\times    y q(y) \lleb(dy),
   \\
   A_3(\delta,M,x) &\eqdef& \int_{\{y, |y| \leq M,   | y' \Gamma_2^{-1} x  | \leq \delta |x| \}}
   \biggl\{   \biggl(  \frac{\pi(x +y) }{ \pi(x)}  -1  \biggr) \1_{R_{x}}(y) + \1_{R_{\infty,x}}(y)  \biggr\}
   \\
   &&\hspace*{33.3mm} {}\times    y q(y) \lleb(dy),
   \\
   A_4(\delta,M,x) &\eqdef& \int_{\{y, |y| \geq M \}}  \biggl\{  \biggl( \frac{\pi(x
        +y) }{ \pi(x)} -1  \biggr) \1_{R_{x}}(y) + \1_{R_{\infty,x}}(y)  \biggr\} y q(y) \lleb(dy).
\end{eqnarray*}
For $x=(x_1,x_2)$ such that $|x_1|-|x_2| \geq 2 M$ and $|y| \leq M$, and $|x_1+y_1| \geq
|x_1| - M \geq |x_2|+M \geq |x_2+y_2|$, it is easily shown that
\begin{eqnarray}\label{eq:borne-pi}
\nonumber  && (1-\alpha) \exp (-0.5 y'\Gamma_2^{-1}  y -    x'\Gamma_2^{-1}  y  )
\\
&&\qquad  \leq \frac{\pi(x +y) }{ \pi(x)}
\\
\nonumber
&& \qquad \leq (1-\alpha)^{-1} \exp (-0.5 y'\Gamma_2^{-1}  y -   x'\Gamma_2^{-1}  y ).
\end{eqnarray}
If $y \in R_{\infty,x} \cap \{ z\dvtx  |x' \Gamma_2^{-1} z | \geq \delta |x|\}$,
then, by \eqref{eq:borne-pi}, $\pi(x + y)/ \pi(x) \leq (1-\alpha)^{-1} \rme^{- \delta
  |x|}$, which implies that $ | A_1(\delta,M,x)  | \leq
(1-\alpha)^{-1} \rme^{- \delta |x|}   \int |y| q(y) \lleb(dy)$. Furthermore,
for any $K$ such that $(1-\alpha)^{-1} \rme^{-\delta K} \leq 1$ and $x$ such
that $ | |x_1|-|x_2|  | \geq 2 M$ and $|x| \geq K$, $R_{\infty,x}
\cap  \{y: |y| \leq M, |x' \Gamma_2^{-1} y| \geq \delta |x|  \}
\subseteq \mathsf{R}_x$.  This property yields to the bound
\begin{eqnarray}
\label{eq2:Borne-A2}
&& \bigg| \frac{\pi( x +y) }{ \pi( x)}  -1  \bigg|  | \1_{R_{ x}}(y) - \1_{R_{\infty,x}}(y)  | \1 {\{y,  |x' \Gamma_2^{-1} y| \geq \delta |x|, |y| \leq M \}}
\nonumber
\\[-8pt]
\\[-8pt]
\nonumber
&&\qquad \leq \1_{R_{x} \setminus R_{\infty,x}}(y) \1 {\{y, |y| \leq M, |x'
  \Gamma_2^{-1} y| \geq \delta |x|\}}.
\end{eqnarray}
Again using \eqref{eq:borne-pi} for $y \in \mathsf{R}_x \cap \{ |y| \leq M \}$,
$(1-\alpha) \rme^{-0.5 a^2 M^2} \rme^{- x' \Gamma_2^{-1} y} \leq \pi(x+y)/\pi(x)
\leq 1$.  On the other hand, for $y \notin R_{\infty,x}$ satisfying $|x'
\Gamma_2^{-1} y | \geq \delta |x|$, we have $x' \Gamma_2^{-1} y \leq -\delta
|x|$, showing that
\begin{eqnarray*}
&& y \in R_{x} \setminus R_{\infty,x} \cap \{z, |z| \leq M, |x' \Gamma_2^{-1} z|
\geq \delta |x|\}
\\
&&\qquad \Longrightarrow\quad (1-\alpha) \rme^{-0.5 a^2 M^2} \rme^{\delta
  K} \leq \pi(x+y)/\pi(x) \leq 1.
\end{eqnarray*}
For fixed $M$, we choose $K$ such that $(1-\alpha) \rme^{-0.5 a^2 M^2} \rme^{
  \delta K} >1$, which implies that the right-hand side in (\ref{eq2:Borne-A2})
is zero and thus $A_2(\delta,M,x)= 0$.  Finally, consider $A_i(\delta,M,x)$,
$i=3,4$. Noting that
\[
 \bigg|   \biggl(  \frac{\pi(x +y) }{ \pi(x)}  -1  \biggr) \1_{R_{x}}(y) + \1_{R_{\infty,x}}(y)  \bigg| \leq 2,
\]
the proof follows from the bounds
\begin{eqnarray}
\label{eq:bound-A3}
|A_3(\delta,M,x)| &\leq& 2 M \int \1 \{ y,  |y' \Gamma_2^{-1} x | \leq \delta |x| \} |y| q(y) \lleb(dy),
\\
\label{eq:bound-A4}
|A_4(\delta,M,x)| &\leq& 2 \int_{|y|\geq M} |y| q(y) \lleb(dy).
\end{eqnarray}
These terms are arbitrarily small for convenient constants $M$ and $\delta$.

\begin{figure}

\includegraphics{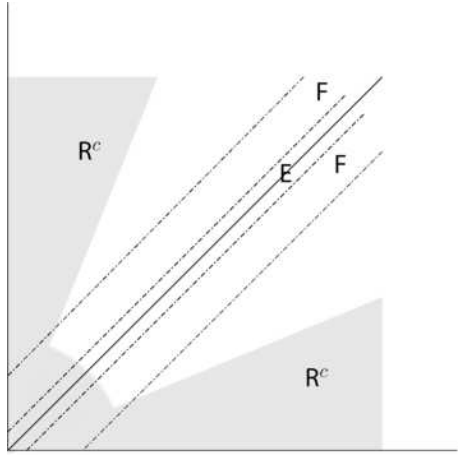}

 \caption{The complement of the zone $\mathsf{R}$ and the strips $\mathsf{E}$ and $\mathsf{F}$.}\label{fig:Zones}
\end{figure}

\subsection[Proof of Proposition 2.10]{Proof of Proposition \textup{\protect\ref{prop:Ex2-cond1PropStab}}}

\subsubsection[Proof of condition (i) of Theorem~1.9]{Proof of condition \textup{(i)} of Theorem~\textup{\protect\ref{theo:stability-fluid-limit-deterministic-nonsmooth}}}
The only difficulty here stems from the irregularity of the ODE for initial
conditions on the diagonals.  Consider the $\beta$-fluid limit
$\QQq{\beta}{u_\star}$ with initial condition $u_\star \eqdef  ( 1/\sqrt{2},
  1/\sqrt{2}  )$ (the other cases can be dealt with similarly).  Set
$v_\star \eqdef  ( 1/\sqrt{2}, -1/\sqrt{2}  )$ and define $V(x) =
|\pscal{v_\star}{x}|$.  Since the increment distribution is assumed to be bounded, there exists a
positive constant $C_q$ such that $|\Phi_1 - \Phi_0| \leq C_q$, $\PP_x$-a.s. for
all $x \in \Xset$.
By Lemma \ref{lem:champ-de-trefle}, we may choose constants
$\gamma \in (0,1)$, $m>0$, $M_0 > C_q$ and $R$ such that\looseness=1
\begin{equation}  \label{eq:InclusionR}
  \mathsf{R} \cap \mathsf{E}^c  \subset \{x \in \Xset,  | \pscal{v_\star}{\Delta(x)} |
\geq m,   \pscal{v_\star}{x} \pscal{v_\star}{\Delta(x)} >0 \},
\end{equation}
where (see Figure~\ref{fig:Zones})
\begin{equation}
\label{eq:definition-R}
\mathsf{E} \eqdef \{x, V(x) \leq M_0 \}\quad\mbox{and}\quad
\mathsf{R}\eqdef  \{ x \in \Xset,|x| \geq R, |\pscal{v_\star}{n(x)}| \leq \gamma  \}.
\end{equation}
For $\delta > 0$, define the stopping time $\kappa(\delta)$ as the infimum of the following three stopping times
\begin{eqnarray}
\label{eq:definition-kappa-1}
\kappa_1(\delta) &\eqdef& \inf \{k \geq 0,  | \pscal{v_\star}{\Phi_k} | \geq 2 \delta  |\Phi_0| \},
\\
\label{eq:definition-kappa-2}
\kappa_2 &\eqdef& \inf \{k \geq 0,  | \Phi_k - \Phi_0  | \geq (1/2) |\Phi_0| \},
\\
\label{eq:definition-kappa-3}
\kappa_3 &\eqdef& \inf \{ k \geq 0, |\Phi_k| < R
\}.
\end{eqnarray}
We will establish the following drift condition: there exist constants $b>0$
and $C$ such that for all $\delta \in (0, \gamma/4)$,
\begin{eqnarray}
\hspace*{10mm}    \PE [V(\Phi_{k+1}) \vert \mathcal{F}_k  ]  &\geq& V(\Phi_k) + m - b\1_{\mathsf{E}}(\Phi_k)\qquad\mbox{on the set $\{k < \kappa(\delta) \}$},     \label{eq:Cond1}
\hspace*{-12pt}\\
   \PE_x  \Biggl[ \sum_{k=0}^{\kappa(\delta) -1} \1_{\mathsf{E}}(\Phi_k) \Biggr]
  &\leq& C,   \label{eq:Step2-Obj}
\end{eqnarray}
with the convention that $\sum_{a}^b =0$ when $a>b$.  We postpone the proof of
\eqref{eq:Cond1} and~\eqref{eq:Step2-Obj} and show how these drift conditions
allow us to obtain condition (i).  On the event
$\{k < \kappa(\delta)\}$, $|\Phi_k| \geq R$, $ (1/2)|\Phi_0| \leq |\Phi_k| \leq
(3/2) |\Phi_0|$ and $ |\pscal{v_\star}{n(\Phi_k)}  |\leq 4 \delta
\leq \gamma$.  Therefore, for all $x \in \Xset$, $\PP_x$-a.s.,
\begin{equation}
\label{eq:BeforeTauInR}
\{ k < \kappa(\delta) \}  \subset \{ \Phi_k \in \mathsf{R}   \}.
\end{equation}
Condition (\ref{eq:Cond1}) yields, for any constant $N>0$,
\[
m   \PE_x[\kappa(\delta) \wedge N] \leq \PE_x\bigl[V\bigl(\Phi_{\kappa(\delta) \wedge
  N}\bigr) \1 \{\kappa(\delta) \geq 1\}\bigr] + b   \PE_x  \Biggl[
  \sum_{k=0}^{\kappa(\delta) \wedge N -1} \mathbh{1}_{\mathsf{E}}(\Phi_k) \Biggr].
\]
The definitions of $\kappa(\delta)$ and $C_q$ imply that
$\PE_x[V(\Phi_{\kappa(\delta) \wedge N}) \1 \{ \kappa(\delta) \geq 1\}] \leq
2 \delta  |x| + C_q $ for all $N$, which, with
(\ref{eq:Step2-Obj}), yields the bound
\begin{equation}
\label{eq:bound-stopping-time}
 m \PE_x [\kappa(\delta)] \leq 2 \delta  |x| + b C  + C_q .
\end{equation}
Let $\{x_n\}$ be a sequence of initial states
such that $\lim_{n \to \infty} x_n = u_\star$ and $\{r_n \}$ be a sequence of
scaling constants, $\lim_{n \to \infty} r_n = +\infty$.  By
Lemma~\ref{lem:ControleEcartLF}, there exists $T_0$ such that
$\QQq{\beta}{u_\star}  \{ \sup_{t \in [0,T_0]} |\eta(t) -
\eta(0)| < 1/4  \}=1$.  Furthermore, we have $1/2 \leq
|x_n| \leq 3/2$ for all $n$ large enough.  Then, by the Portmanteau
theorem,
\begin{eqnarray*}
  && \QQq{\beta}{u_\star} \{ \eta, \eta([0,T_0]) \cap \openset = \varnothing  \}
  \\
  &&\qquad  = \lim_{\delta \downarrow 0^+} \QQq{\beta}{u_\star}
  \biggl\{ \eta, \sup_{t \in [0,T_0]} |\eta(t) - \eta(0)| < 1/4, \sup_{t \in [0,T_0]}  |\pscal{v_\star}{\eta(t)}  |  < \delta  \biggr\}
  \\
  &&\qquad  \leq \lim_{\delta \downarrow 0^+} \liminf_{n \to \infty} \PP_{r_n x_n}
  \biggl\{ \sup_{0 \leq k \leq 2 T_0  |\Phi_0|/3 }  | \Phi_k - \Phi_0  | < (1/2) |\Phi_0|,
\\
&&\hspace*{46.5mm}    \sup_{0 \leq k \leq 2 T_0  |\Phi_0| /3}   | \pscal{v_\star}{\Phi_k}  | < 2 \delta  |\Phi_0| \biggr\}
  \\
  &&\qquad \leq \lim_{\delta \downarrow 0^+} \liminf_{n \to \infty} \PP_{r_n x_n}
   \bigl( \kappa(\delta) \geq 2 T_0 | \Phi_0|/3  \bigr) = 0,
\end{eqnarray*}
where the last equality stems from \eqref{eq:bound-stopping-time}. This proves
Theorem~\ref{theo:stability-fluid-limit-deterministic-nonsmooth}(i).

We now prove (\ref{eq:Cond1}). Since $\PE [ \Phi_{k+1}\vert \mathcal{F}_k
 ] = \Phi_k + \Delta(\Phi_k)$, Jensen's inequality implies that $\PE_x
 [V(\Phi_{k+1}) \vert \mathcal{F}_k  ] \geq  |
  \pscal{v_\star}{\Phi_k+\Delta(\Phi_k)}  | $.  Furthermore, by
\eqref{eq:InclusionR} and \eqref{eq:BeforeTauInR}, $\{ k <\kappa(\delta),
\Phi_k \in \mathsf{E}^c \} \subset \{\Phi_k \in \mathsf{R} \cap \mathsf{E}^c
\}$, which implies that $ | \pscal{v_\star}{\Phi_k+\Delta(\Phi_k)}  | - |
\pscal{v_\star}{\Phi_k} |=  | \pscal{v_\star}{\Delta(\Phi_k)}  | \geq
m$ since, on $\mathsf{R}\cap \mathsf{E}^c$, $\pscal{v_\star}{x}$ and
$\pscal{v_\star}{\Delta(x)}$ have the same sign and
$\pscal{v_\star}{\Delta(x)}$ is lower bounded.  On the set $\{ k<
\kappa(\delta), \Phi_k \in \mathsf{E} \}$, we write $V(\Phi_{k+1}) \geq
V(\Phi_k) - C_q$ so that $\PE  [V(\Phi_{k+1}) \vert \mathcal{F}_k  ]
\geq V(\Phi_k) +m - (C_q+m)$.  This concludes the proof of (\ref{eq:Cond1}).

Finally, we prove (\ref{eq:Step2-Obj}).  For $A \in \Xsigma$, we denote by
$\sigma_A \eqdef \inf\{k \geq 0, \Phi_k \in A\}$ the first hitting time on $A$.
For notational simplicity, we write $\kappa$ instead of $\kappa(\delta)$.
Define recursively $\sigma^{(1)} \eqdef \sigma_{\mathsf{E}\cap \mathsf{R}}$
and, for all $k \geq 2$, $\sigma^{(k)} \eqdef \sigma^{(k-1)} + \tau \circ
\theta^{\sigma^{(k-1)}} + \sigma^{(1)} \circ \theta^{\tau \circ
  \theta^{\sigma^{(k-1)}} + \sigma^{(k-1)}}$, where $\tau \eqdef \kappa \wedge
k_\star$, $k_\star$ being an integer whose value will be specified later.
With this notation,
\begin{equation}
  \label{eq:UpperBoundNbrVisitE}
  \PE_x  \Biggl[\sum_{k=0}^{\kappa -1} \1_{\mathsf{E}}(\Phi_k)  \Biggr] \leq k_\star \sum_{q \geq 1} \PP_x  \bigl( \sigma^{(q)} < \kappa  \bigr).
\end{equation}
Furthermore, for all $q \geq 2$, the strong Markov property yields the bound
\[
\PP_x  \bigl( \sigma^{(q)} < \kappa  \bigr) \leq \PP_x  \bigl(\sigma^{(q-1)} <
  \kappa  \bigr) \sup_{y \in \mathsf{E} \cap \mathsf{R}} \PP_y  \bigl(\tau +
  \sigma^{(1)} \circ \theta^{\tau} < \kappa  \bigr).
\]
Therefore, by (\ref{eq:UpperBoundNbrVisitE}), (\ref{eq:Step2-Obj}) holds, provided that
$\sup_{x \in \mathsf{E} \cap \mathsf{R}} \PP_x  (\tau + \sigma^{(1)} \circ \theta^{\tau} < \kappa  ) <1$.  For all $x \in \mathsf{E} \cap \mathsf{R}$,
it is easily seen that
\begin{eqnarray}
&&  \PP_x  \bigl( \tau + \sigma^{(1)} \circ \theta^{\tau} < \kappa  \bigr)
\nonumber
\\[-8pt]
\\[-8pt]
\nonumber
&&\qquad =
  \PP_x ( \tau < \kappa) -
  \PE_x  \bigl(\1 \{ \tau < \kappa \} \1 \{ \Phi_\tau \in \mathsf{E}^c \cap \mathsf{R} \} \PP_{\Phi_\tau}  \bigl[\kappa \leq \sigma^{(1)}  \bigr]  \bigr)
  \\
&&\qquad \leq  1 - \inf_{x \in \mathsf{E}^c \cap \mathsf{R}} \PP_x  \bigl( \kappa \leq
    \sigma^{(1)}  \bigr)  \{ \PP_x  ( \tau = \kappa  ) + \PP_x
     ( \tau= k_\star, \Phi_{k_\star} \in \mathsf{E}^c \cap \mathsf{R}
     )  \},\hspace*{-12pt}
\end{eqnarray}
showing that the conditions
\begin{eqnarray}
\inf_{x \in \mathsf{E} \cap \mathsf{R}} \PP_x  ( \{ \tau < k_\star \} \cup \{\tau= k_\star, \Phi_{k_\star} \in \mathsf{E}^c \cap \mathsf{R} \} )  &>& 0,   \label{eq:Cond2}
\\
\inf_{x \in \mathsf{E}^c \cap \mathsf{R}} \PP_x  \bigl( \kappa \leq \sigma^{(1)}  \bigr) &>& 0  \label{eq:Cond4}
\end{eqnarray}
imply \eqref{eq:Step2-Obj}. We first prove \eqref{eq:Cond2}. Choose  $\tilde \gamma \in (\gamma,1)$ such that the four half-planes $\{z,
\pscal{z}{\Gamma_i^{-1}u_{\star,\tilde \gamma}^\pm} < 0\}$ ($i=1,2$) have a nonempty intersection, where
$u_{\star,\tilde{\gamma}}^{-}$ and $u_{\star, \tilde{\gamma}}^+$ are the unit vectors defining the
edges of the cone $\mathsf{C}_{\tilde \gamma} \eqdef \{z \in \Xset, |\pscal{v_\star}{n(z)} \leq \tilde \gamma \}$.
Define
\begin{equation}
\label{eq:definition-W}
\mathsf{W} \eqdef \{z, 0 \leq |z| \leq C_q,  \pscal{z}{\Gamma_i^{-1}u^{\pm}_{\star,\tilde \gamma}} \leq 0,   i=1,2 \}.
\end{equation}
Since any vector $y$ in the cone $\mathsf{C}_{\tilde \gamma}$ can be written as
a linear combination of the vectors $u_{\star,\tilde \gamma}^-$ and $u_{\star, \tilde \gamma}^+$ with positive
weights,
for any $y \in \mathsf{C}_{\tilde \gamma}$ and $z \in \mathsf{W}$, $\pscal{z}{\Gamma_i^{-1} y} \leq 0$, $i=1,2$, which implies that
\begin{eqnarray*}
&& \pscal{z}{\nabla \pi(y)}
\\
&&\qquad = -\alpha \pscal{z}{\Gamma_1^{-1}y} \exp( - 0.5 y' \Gamma_1^{-1} y) -  (1-\alpha) \pscal{z}{\Gamma_2^{-1}y} \exp( - 0.5 y' \Gamma_2^{-1} y)
\\
&&\qquad \geq 0.
\end{eqnarray*}
By choosing $R$ large enough [see (\ref{eq:definition-R})], we can
assume,  without loss of generality, that for all $x \in \mathsf{R}$ and $z \in
\mathsf{W}$, $x + tz \in \mathsf{C}_{\tilde \gamma}$ for all $t \in (0,1)$.
Thus, $\pi(x+z) = \pi(x) + \int_0^1 \pscal{\nabla \pi(x+tz)}{z} \,dt \geq 0$ and
we have $\pi(x+z) \geq \pi(x)$, showing that $\mathsf{W} \subset \mathsf{A}_x$.
Finally, we write $\mathsf{W}$ as the union of two disjoint sets
$\mathsf{W}^-$, $\mathsf{W}^+$, where $\mathsf{W}^+ \eqdef \{z \in \mathsf{W},
\pscal{v_\star}{z} \geq 0 \}$.  Since, for $x \in \mathsf{R}$, $\mathsf{W}
\subset \mathsf{A}_x$, for any $0 \leq c \leq C_q$, we have
\begin{eqnarray*}
&& \inf_{x \in \mathsf{R}, \pscal{v_\star}{x} \geq 0} \PP_x \bigl(
  |\pscal{v_\star}{\Phi_1}| \geq |\pscal{v_\star}{\Phi_0}| + c  \bigr)
\\
&&\qquad  \geq
\int_{\mathsf{W}^+} \1 {\{y, |\pscal{v_\star}{y}| \geq c\}} q(y) \lleb(dy) >0.
\end{eqnarray*}
An analogous lower bound holds for all $x \in \mathsf{R}$ such that
$\pscal{v_\star}{x} \leq 0$.  These inequalities, combined with repeated
applications of the Markov property, yield \eqref{eq:Cond2}, by choosing
$k_\star$ such that $k_\star c \geq M_0$.

We now prove \eqref{eq:Cond4}.
Let $M_1 > M_0$ and set $\mathsf{F} \eqdef \{x, V(x) \leq M_1 \}$.
By Lemma~\ref{lem:borne-M-infty}, we may choose  $J\geq 1$ and then $M_1 >M_0$ large enough so that, for all $x \in \Xset$,
\begin{eqnarray}\label{eq:definition-J}
\PP_x  \Biggl( \sup_{j \geq J} j^{-1}  \Bigg| \sum_{l=1}^j \epsilon_l  \Bigg| \geq
  m  \Biggr)< 1/2,
\nonumber
\\[-8pt]
\\[-8pt]
\nonumber
\PP_x  \Biggl( \sup_{j \leq J}  \Bigg| \sum_{l=1}^j
    \epsilon_l  \Bigg| \geq M_1 - M_0  \Biggr)< 1/2.
\end{eqnarray}
It is easily seen that, using the strong Markov property,
\[
\inf_{x \in \mathsf{E}^c \cap \mathsf{R}} \PP_x  (\kappa \leq \sigma_{\mathsf{E} \cap \mathsf{R}}  ) \geq
\inf_{x \in \mathsf{E}^c \cap \mathsf{R}} \PP_x  (\sigma_{\mathsf{F}^c \cap \mathsf{R}} < \sigma_{\mathsf{E} \cap \mathsf{R}}  )
\inf_{x \in \mathsf{F}^c \cap \mathsf{R}} \PP_x  (\kappa \leq \sigma_{\mathsf{E}\cap \mathsf{R}}  ).
\]
The first term of the right-hand side of the previous relation can be shown to be positive, using arguments
which are similar to those used in the proof of \eqref{eq:Cond2}.  We write $
\pscal{v_\star}{\Phi_{k}}= \pscal{v_\star}{\Phi_{0}} + \sum_{l=1}^{k}
\pscal{v_\star}{\Delta(\Phi_{l-1})} +\sum_{l=1}^{k}
\pscal{v_\star}{\epsilon_l}$. Let $x \in \mathsf{F}^c \cap \mathsf{R}$.
$\PP_x$-a.s., since $|\Phi_l - \Phi_{l-1}| \leq C_q \leq M_0$, on the event $\{ 1 \leq
k \leq \sigma_{\mathsf{E} \cap \mathsf{R}} < \kappa \}$,
$|\pscal{v_\star}{\Phi_k}| \geq M_0$, $\pscal{v_\star}{\Phi_{0}}
\pscal{v_\star}{\Phi_j} > 0$ and $\pscal{v_\star}{\Phi_0}
\pscal{v_\star}{\Delta(\Phi_j)} >0$ for all $0 \leq j < k$, which
implies that
\begin{eqnarray*}
 | \pscal{v_\star}{\Phi_{k}}  |
&\geq&  | \pscal{v_\star}{\Phi_{0}}  | + \sum_{l=1}^{k}  | \pscal{v_\star}{\Delta(\Phi_{l-1})}  | -  \Bigg|\sum_{l=1}^{k} \pscal{v_\star}{\epsilon_l}  \Bigg|
\\
&\geq& M_1 + k m -  \Bigg|\sum_{l=1}^{k} \pscal{v_\star}{\epsilon_l}  \Bigg|.
\end{eqnarray*}
Thus, for all $x \in \mathsf{F}^c \cap \mathsf{R}$, using the definition
\eqref{eq:definition-J} of $J$ and $M_1$, we have
\begin{eqnarray*}
  \PP_x \{ J \leq  \sigma_{\mathsf{E} \cap \mathsf{R}} < \kappa \}
  &\leq & \sup_{x \in \Xset} \PP_x  \Biggl\{ \sup_{j \geq J} j^{-1}  \Bigg|\sum_{l=1}^{j} \pscal{v_\star}{\epsilon_l} \Bigg| \geq m  \Biggr\} < 1/2,
  \\
  \PP_x \{ \sigma_{\mathsf{E} \cap \mathsf{R}} < \kappa \wedge J \} &\leq & \sup_{x \in \Xset}
  \PP_x  \Biggl\{ \sup_{j \leq J}  \Bigg|\sum_{l=1}^{j} \epsilon_l  \Bigg| \geq (M_1-M_0)  \Biggr\} < 1/2,
\end{eqnarray*}
which proves $\inf_{x \in \mathsf{F}^c \cap \mathsf{R}} \PP_x  (\kappa \leq \sigma_{\mathsf{E}\cap \mathsf{R}}  )>0$ and therefore \eqref{eq:Cond4}.

\subsubsection[Proof of  B4 and   the conditions (ii)--(iii) of Theorem~1.9]{Proof of  \textup{B4} and   the conditions \textup{(ii)--(iii)}
of Theorem~\textup{\protect\ref{theo:stability-fluid-limit-deterministic-nonsmooth}}}
Assume that $x \in \mathsf{C} \eqdef \{x, 0<|x_2|<x_1\}$ (the three other cases
are similar). By Lemma~\ref{lem:champ-de-trefle}, $h(x) = - c_q n(\Gamma_2^{-1}
x)$ for all $x \in \mathsf{C}$, which is locally Lipschitz. Hence, there exists
a unique maximal solution $\mu(\cdot;x)$ on $[0,T_x]$ satisfying
$\mu(0;x)=x$ and $\mu(t;x) \in \mathsf{C}$ for all $t \leq
T_x$, showing B4.  Since, for $t \in [0,T_x)$, $d/dt
|\mu(t;x)|^2 = 2 |\mu(t;x)| \pscal{n(\mu(t;x))}{h
  \circ \mu(t;x)} < -2 c_q |a|^{-1} |\mu(t;x)|$, the norm
of the ODE solution is bounded by $|\mu(t;x)| \leq (|x|- c_q |a|^{-1}
t)_+$ for all $0 \leq t \leq |x| |a| c_q^{-1}$, which implies condition
(ii), provided that $T_x \geq T_K$ for all $x \in \mathsf{C}
\cap \mathsf{B}(0,K)$. This result follows from the fact that the boundaries of
$\mathsf{C}$ are repulsive: consider the relative neighborhood in $\mathsf{C}$,
$\mathsf{V} \eqdef \mathsf{V}_1 \cup \mathsf{V}_2$, of the boundaries where
$\mathsf{V}_1 \eqdef \{ x \dvtx  x_1 > 0, \pscal{v_\star}{x} > 0,
\pscal{x}{\Gamma_2^{-1} v_\star} < 0\}$ and $\mathsf{V}_2 \eqdef \{x \dvtx x_1 > 0,
\pscal{u_\star}{x} > 0, \pscal{x}{ \Gamma_2^{-1} u_\star} < 0 \}$. Assume that
there exists $s \in [0, T_x]$ such that $\mu(s;x) \in \mathsf{V}_1$
(the other case can be handled similarly).  Since $t \mapsto
\mu(t;x)$ is continuous and $\mathsf{V}_1$ is a relative open subset
of $\mathsf{C}$, there exists $\delta$ such that for all $0 \leq t \leq
\delta$, $\mu(s+t;x) \in \mathsf{V}_1$. This implies that for all $0
\leq t \leq \delta$,
\begin{eqnarray*}
&& \pscal{v_\star}{ \mu(s+t;x)} - \pscal{v_\star}{\mu(s;x)}
\\
&&\qquad =
-c_q \int_0^t |\Gamma_2^{-1} \mu(s+u;x)|^{-1}
\pscal{v_\star}{ \Gamma_2^{-1} \mu(s+u;x)} \,du >0,
\end{eqnarray*}
showing that, in $\mathsf{V}_1$, the distance to the boundary always increases.
The properties above also imply condition (iii) of
Theorem~\ref{theo:stability-fluid-limit-deterministic-nonsmooth}.

\begin{appendix} 
\section*{Appendix: Technical lemmas}

\begin{lem}
\label{lem:borne-M-infty}
Let $\{ \varepsilon_k \}_{k \geq 1}$ be an $L^p$-martingale difference sequence
adapted to the filtration $\{ \mathcal{F}_k \}_{k \geq 0}$.  For any $p > 1$,
there exists a constant $C$ (depending only on $p$) such that
\begin{eqnarray}
\label{eq:doob-maximal-inequality}
\PE  \Biggl[ \sup_{1 \leq l \leq n} \Bigg| \sum_{k=1}^l \varepsilon_k  \Bigg|^p  \Biggr]
&\leq& C \sup_{k \geq 1} \PE[ |\varepsilon_k|^p ] n^{1 \vee p/2},
\\
\label{eq:birnbaum-marshall-inequality}
\PP  \Biggl[ \sup_{n \leq l} l^{-1} \Bigg| \sum_{k=1}^l \varepsilon_k  \Bigg|
\geq M  \Biggr] &\leq& C \sup_{k \geq 1} \PE[ |\varepsilon_k|^p ] M^{-p} n^{-p+1\vee p/2}.
\end{eqnarray}
\end{lem}

\begin{pf}
For $p>1$, applying in sequence the Doob maximal inequality, by the
Burkholder inequality for $L^p$-martingales, there exists a constant $C_p$ such that
\[
\PE  \Biggl[  \sup_{1 \leq l \leq n}  \Bigg| \sum_{k=1}^l \varepsilon_k  \Bigg|^p  \Biggr]
\leq C_p \PE  \Biggl[  \Bigg| \sum_{k=1}^{n}  |\varepsilon_k  |^2  \Bigg|^{p/2}  \Biggr].
\]
Equation \eqref{eq:doob-maximal-inequality} follows from the Minkovski inequality for $p\geq 2$,
\begin{equation}
\label{eq:borne-M-infty-p>2}
\PE  \Biggl[  \sup_{1 \leq l \leq n} \Bigg| \sum_{k=1}^l \varepsilon_k  \Bigg|^p  \Biggr]
\leq  C_p  \sup_{k \geq 1} \PE[ |\varepsilon_k|^p ] n^{p/2},
\end{equation}
and the subadditivity inequality for $ 1 < p \leq 2$,
\begin{equation}
\label{eq:borne-M-infty-1<p<2}
\PE  \Biggl[ \sup_{1 \leq l \leq n}  \Bigg| \sum_{k=1}^l \varepsilon_k  \Bigg|^p   \Biggr]
\leq  C_p  \sup_{k \geq 1} \PE[ |\varepsilon_k|^p ]  n.
\end{equation}
Equation~\eqref{eq:birnbaum-marshall-inequality} follows from \citet{birnbaummarshall1961}, Theorem~1.
\end{pf}

\begin{lem}
\label{lem:elem-1}
Let $X,Y$ be two nonnegative random variables. Then,
for any $p \geq 1$, there exists a constant $C_p$ (depending only on $p$) such that, for any $M > 0$,
\[
\PE[ (X+Y)^p \1 \{ X+Y > M \}] \leq C_p  ( \PE  [ X^p \1 \{ X \geq M/2 \}  ] + \PE[ Y^p ]  ).
\]
\end{lem}

\begin{pf}
Note that $\1 \{ X + Y \geq M \} \leq \1 \{ X \geq M/2 \} + \1 \{ X \leq M/2 \} \1 \{ Y \geq M/2 \}$. Therefore,
\begin{eqnarray*}
\PE( X^p \1 \{ X+Y \geq M \}) &\leq& \PE( X^p \1 \{ X \geq M/2 \}) + (M/2)^p \PP( Y \geq M/2)
\\
&\leq& \PE( X^p \1 \{ X \geq M/2 \}) + \PE(Y^p).
\end{eqnarray*}
The proof then follows from the fact that $(X+Y)^p \leq 2^{p-1} (X^p+Y^p)$.
\end{pf}

\begin{lem}
\label{lem:elem-2} Let $X$ be a nonnegative random variable. For any $p
\geq 0$, $a > 1$ and $M$, we have
\[
\PE[ X^p \1 \{ X \geq M \}] \leq M^{-(a-1)p} \PE[ X^{ap}].
\]
\end{lem}

\begin{pf}
\begin{eqnarray*}
\PE[ X^p \1 \{ X \geq M \}] &\leq&  ( \PE[X^{ap}]  )^{1/a}
( \PP[ X \geq M]  )^{(a-1)/a}
\\
&\leq&  ( \PE[X^{ap}]  )^{1/a}  ( M^{-ap} \PE[X^{ap}]  )^{(a-1)/a}.
\end{eqnarray*}\upqed
\end{pf}

\end{appendix}

\printaddresses

\end{document}